\documentclass[11pt]{article}
\usepackage{amssymb}
\usepackage{amsmath}
\usepackage{amsthm}
\newcommand{\cA}{{\cal A}}
\newcommand{\cB}{{\cal B}}
\newcommand{\cC}{{\cal C}}

\newcommand{\cH}{{\cal H}}
\newcommand{\cE}{{\cal E}}

\newcommand{\cI}{{\cal I}}
\newcommand{\cO}{{\cal O}}
\newcommand{\cL}{{\cal L}}
\newcommand{\cM}{{\cal M}}
\newcommand{\cN}{{\cal N}}
\newcommand{\cF}{{\cal F}}

\newcommand{\cV}{{\cal V}}
\newcommand{\cW}{{\cal W}}
\newcommand{\cX}{{\cal X}}
\newcommand{\cY}{{\cal Y}}

\renewcommand{\AA}{{\mathbb A}}

\newcommand{\ZZ}{{\mathbb Z}}

\newcommand{\gp}{\mathfrak{p}}
\newcommand{\gq}{\mathfrak{q}}

\newcommand{\gt}{\mathfrak{t}}

\newcommand{\gF}{\mathfrak{F}}

\newcommand{\on}{\operatorname}

\newcommand{\Rep}{{\on{Rep}}}

\newcommand{\Qlb}{\mathbb{\bar Q}_\ell}
\newcommand{\Gm}{\mathbb{G}_m}

\newcommand{\A}{\mathbb{A}}

\newcommand{\toup}[1]{\stackrel{#1}{\to}}
\newcommand{\hook}[1]{\stackrel{#1}{\hookrightarrow}}

\newcommand{\getsup}[1]{\stackrel{#1}{\gets}}
\newcommand{\Sp}{\on{\mathbb{S}p}}

\newcommand{\GSp}{\on{G\mathbb{S}p}}

\newcommand{\Hom}{\on{Hom}}

\newcommand{\Sym}{\on{Sym}}

\newcommand{\Ker}{\on{Ker}}

\newcommand{\Aut}{\on{Aut}}

\newcommand{\Pic}{\on{Pic}}

\newcommand{\Bun}{\on{Bun}}
\newcommand{\Bunb}{\on{\overline{Bun}} }
\newcommand{\Bunt}{\on{\widetilde\Bun}}

\newcommand{\Spec}{\on{Spec}}

\newcommand{\HOM}{{{\cal H}om}}
\newcommand{\END}{{{\cal E}nd}}
\newcommand{\Gr}{\on{Gr}}
\newcommand{\GR}{\on{{\mathcal{G}}r}}
\newcommand{\Grb}{\overline{\Gr}}

\newcommand{\GL}{\on{GL}}

\newcommand{\Gal}{\on{Gal}}

\newcommand{\pr}{\on{pr}}
\newcommand{\id}{\on{id}}
\newcommand{\tr}{\on{tr}}
\newcommand{\QED}{$\square$} 
\newcommand{\Fq}{\mathbb{F}_q}  
\newcommand{\Fp}{\mathbb{F}_p}  
\newcommand{\iso}{{\widetilde\to}}

\newcommand{\comp}{\circ}

\renewcommand{\H}{{\on{H}}}   

\newcommand{\DD}{\mathbb{D}}  
\newcommand{\D}{\on{D}}       

\newcommand{\overl}[1]{\overline{#1}}
\newcommand{\select}[1]{{\it{#1}}}

\newcommand{\ov}[1]{\overline{#1}}

\renewcommand{\div}{\on{div}}

\renewcommand{\P}{{\on{P}}}

\newcommand{\<}{\langle}
\renewcommand{\>}{\rangle}

\newcommand{\ev}{\on{ev}}

\newcommand{\Sph}{\on{Sph}}
\newcommand{\Res}{\on{Res}}

\newcommand{\tboxtimes}{\tilde\boxtimes}

\newcommand{\act}{\on{act}}

\newcommand{\RCov}{\on{RCov}}

\newcommand{\diag}{\on{diag}}

\newcommand{\Wald}{\on{{\cal W}ald}}

\newcommand{\Isom}{\on{Isom}}

\newcommand{\BM}{\on{BM}}
\newcommand{\Orb}{\on{Orb}}

\newtheorem{Lm}{Lemma}
\newtheorem{Th}{Theorem}
\newtheorem{Pp}{Proposition}
\newtheorem{Cor}{Corolary}

\theoremstyle{remark}
\newtheorem{Rem}{Remark}

\theoremstyle{definition}
\newtheorem{Def}{Definition}

\newenvironment{Prf}{\par\noindent {\it Proof }}{\QED}
\newcommand{\Step}[1]{\par\noindent{\bf Step {#1}}.}

\begin{document}

\author{Sergey Lysenko}
\title{Geometric Bessel models for $\GSp_4$ and multiplicity one}
\date{}
\maketitle
\begin{abstract}
\noindent{\scshape Abstract}\hskip 0.8 em  I this paper, which is a sequel to \cite{Ly2}, we study Bessel models of representations of $\GSp_4$ over a local non archimedian field in the framework of the geometric Langlands program. 
 The Bessel module over the nonramified  Hecke algebra of $\GSp_4$ admits a geometric counterpart, the Bessel category of perverse sheaves on some ind-algebraic stack.  We use it to prove a geometric version of multiplicity one for Bessel models. It implies a geometric Casselman-Shalika type formula for these models. 
 The strategy of the proof is the same as in a paper of Frenkel, Gaitsgory and Vilonen \cite{FGV}. 
 We also propose a geometric framework unifying Whittaker, Waldspurger and Bessel models.
\end{abstract} 

\medskip\smallskip

{\centerline{\scshape 0. Introduction}}

\medskip\noindent
0.1  In this paper, which is a sequel to \cite{Ly2}, we study Bessel models of representations of $\GSp_4$ in the framework of the geometric Langlands program. These models introduced by Novodvorsky and Piatetski-Shapiro, satisfy the following multiplicity one property (\cite{NP}).

 Set $k=\Fq$ and $\cO=k[[t]]\subset F=k((t))$. Let $\tilde F$ be an \'etale $F$-algebra with $\dim_F(\tilde F)=2$ such that $k$ is algebraically closed in $\tilde F$. Write $\tilde\cO$ for the integral closure of $\cO$ in $\tilde F$. We have two cases:
\begin{itemize} 
\item $\tilde F\,\iso\, k((t^{\frac{1}{2}}))$ (nonsplit case)
\item $\tilde F\,\iso\, F\oplus F$ (split case) 
\end{itemize}
 Write $L$ for $\tilde\cO$ viewed as $\cO$-module, it is equipped with a quadratic form $s:\Sym^2 L\to\cO$ given by the determinant. Write $\Omega_{\cO}$ for the completed module of relative differentials of $\cO$ over $k$.  

 Set $\cM=L\oplus (L^*\otimes\Omega_{\cO}^{-1})$. This $\cO$-module is equipped with a symplectic form $\wedge^2 \cM\to L\otimes L^*\otimes\Omega_{\cO}^{-1}\to\Omega_{\cO}^{-1}$.  
Set $G=GSp(\cM)$, this is a group scheme over $\Spec\cO$. Write $P\subset G$ for the Siegel parabolic subgroup preserving the lagrangian submodule $L$. Its unipotent radical $U$ has a distinguished character 
$$
\ev: U\,\iso\,\Omega_{\cO}\otimes\Sym^2 L\toup{s} \Omega_{\cO}
$$ 
(here we view $\Omega_{\cO}$ as a commutative group scheme over $\Spec\cO$). Set 
$$
\tilde R=\{p\in P\mid ev(pup^{-1})=\ev(u)\;\mbox{for}\; u\in U\}
$$ 
View $\GL(L)$ as a group scheme over $\Spec\cO$ and $\tilde\cO^*$ as its closed subgroup. Write $\alpha$ for the composition $\tilde\cO^*\hook{}\GL(L)\toup{\det}\cO^*$. Fix a section
$\tilde\cO^*\hook{}\tilde R$ given by $g\mapsto (g, \alpha(g)(g^*)^{-1})$. Then $R=\tilde\cO^*U\subset\tilde R$ is a closed subgroup, and the map $R\toup{\xi} \Omega_{\cO}\times \tilde\cO^*$ sending $tu$ to $(ev(u), t)$ is a homomorphism of group schemes over $\Spec\cO$.

 Let $\ell$ be a prime invertible in $k$. Fix a character $\chi:\tilde F^*/\tilde\cO^*\to\Qlb^*$
and a nontrivial additive character $\psi:k\to\Qlb^*$. Write $\tau$ for the composition
$$
R(F)\toup{\xi} \Omega_F\times \tilde F^*\toup{\Res\times\pr} k\times\tilde F^*/\tilde\cO^*\toup{\psi\times\chi} \Qlb^*
$$
The Bessel module is the vector space
\begin{multline*}
\BM_{\tau}=\{f:G(F)/G(\cO)\to\Qlb\mid f(rg)=\tau(r)f(g)\;\mbox{for}\; r\in R(F),\\ 
f\;\mbox{is of compact support modulo}\; R(F)\}
\end{multline*}
Let $\chi_c:F^*/\cO^*\to\Qlb^*$ denote the restriction of $\chi$. The Hecke algebra
\begin{multline*}
\H_{\chi_c}=\{h: G(\cO)\backslash G(F)/G(\cO)\to\Qlb\mid h(zg)=\chi_c(z)h(g)\;\mbox{for}\; z\in F^*, \\
h\;\mbox{is of compact support modulo} \; F^*\}
\end{multline*}
acts on $\BM_{\tau}$ by convolutions. Then $\BM_{\tau}$ is \select{a free module of rank one} over $\H_{\chi_c}$. In this paper we prove a geometric version of this result.  

 Remind that the affine grassmanian $\Gr_G=G(F)/G(\cO)$ can be viewed as an ind-scheme over $k$. According to `fonctions-faisceaux' philosophy, the space $\BM_{\tau}$ should have a geometric counterpart. A natural candidate for that would be the category of $\ell$-adic perverse sheaves on $\Gr_G$ that change under the action of $R(F)$ by $\tau$. However, the $R(F)$-orbits on $\Gr_G$ are infinite-dimensional, and this naive definition does not make sense. 
 
 The same difficulty appears when one tries to define Whittaker categories for any reductive group. In \cite{FGV} Frenkel, Gaitsgory and Vilonen have overcomed this by replacing the corresponding local statement by its globalization, which admits a geometric counterpart leading to a definition of Whittaker categories with expected properties. We follow the strategy of \select{loc.cit.} replacing the above local statement by a global one, which we further geometrize.
    
\medskip\noindent    
0.2  Fix a smooth projective absolutely irreducible curve $X$ over $k$. Let 
$\pi: \tilde X\to X$ be a two-sheeted covering ramified at some effective divisor $D_{\pi}$ of $X$ (we assume $\tilde X$ smooth over $k$). The vector bundle $L=\pi_*\cO_{\tilde X}$ is equipped with a quadratic form $s:\Sym^2 L\to\cO_X$. 

 Write $\Omega$ for the canonical line bundle on $X$. Set $\cM=L\oplus (L^*\otimes\Omega^{-1})$, it is equipped with a symplectic form 
$$
\wedge^2\cM\to L\otimes L^*\otimes\Omega^{-1}  \to\Omega^{-1}
$$ 
Let $G$ be the group scheme (over $X$) of automorphisms of $\cM$ preserving this symplectic form up to a multiple. Let $P\subset G$ denote the Siegel parabolic subgroup preserving $L$, $U\subset P$ its unipotent radical. Then $U$ is equipped with a homomorphism of group schemes over $X$
$$
\ev: U\,\iso\,\Omega\otimes\Sym^2 L\toup{s}\Omega
$$
Let $T$ be the functor sending a $X$-scheme $S$ to the group $\H^0(\tilde X\times_X S, \cO^*)$. Then $T$ is a group scheme over $X$, a subgroup of $\GL(L)$. Write $\alpha$  for the composition $T\hook{}\GL(L)\toup{\det}\Gm$. Set 
$$
\tilde R=\{p\in P\mid \ev(pup^{-1})=\ev(u)\;\,\mbox{for all}\; u\in U\}
$$ 
Fix a section $T\hook{}\tilde R$ given by $g\mapsto (g, \alpha(g)(g^*)^{-1})$. Then $R=TU\subset\tilde R$ is a closed subgroup, and the map $R\toup{\xi} \Omega\times T$ sending $tu$ to $(ev(u), t)$ is a homomorphism of group schemes over $X$.
 
  Let $F=k(X)$, $\AA$ be the adele ring of $F$ and $\cO\subset \AA$ the entire adeles.  Write $F_x$ for the completion of $F$ at $x\in X$ and $\cO_x\subset F_x$ for its ring of integers.  
Fix a nonramified character $\chi: T(F)\backslash T(\AA)/T(\cO)\to \Qlb^*$. Let $\tau$ be the composition
$$
R(\AA)\toup{\xi} \Omega(\AA)\times T(\AA)\,\toup{r\times\chi} \,\Qlb^*,  
$$ 
where $r: \Omega(\AA)\to\Qlb^*$ is given by 
$$
r(\omega_x)=\psi(\sum_{x\in X} \tr_{k(x)/k}\Res \omega_x)
$$
 
 Fix $x\in X(k)$. Let $Y$ denote the restricted product 
$G(F_x)/G(\cO_x)\times \mathop{\prod'}\limits_{y\ne x} R(F_y)/R(\cO_y)$. Let $\cY(k)$ be the quotient of $Y$ by the diagonal action of $R(F)$. Set
\begin{multline*}
\BM_{X,\tau}=\{f: Y\to\Qlb\mid f(rg)=\tau(r)f(g)\;\mbox{for}\; r\in R(\AA),\\ 
 f \;\mbox{is of compact support modulo}\; R(\AA)\}
\end{multline*} 
 View elements of $\BM_{X,\tau}$ as functions on $\cY(k)$. 
Let $\chi_c: F_x^*/\cO_x^*\to\Qlb^*$ be the restriction of $\chi$. As in 0.1, the Hecke algebra $\H_{\chi_c}$ of the pair $(G(F_x), G(\cO_x))$ acts on $\BM_{X,\tau}$ by convolutions. The restriction under
$$
G(F_x)/G(\cO_x)\hook{} Y
$$ 
yields an isomorphism of $\H_{\chi_c}$-modules $\BM_{X,\tau}\to \BM_{\tau}$. 
  
 We introduce an ind-algebraic stack $_{x,\infty}\Bunb_{R_{\pi}}$ whose set of $k$-points contains $\cY(k)$. We define the Bessel category $\P^{\cL}(_{x,\infty}\Bunb_{R_{\pi}})$, a category of perverse sheaves on $_{x,\infty}\Bunb_{R_{\pi}}$ with some equivariance property. This is a geometric version of $\BM_{X,\tau}$. 
  
  Let $\Sph(\Gr_G)$ denote the category of $G(\cO_x)$-equivariant perverse sheaves on the affine grassmanian $G(F_x)/G(\cO_x)$. By \cite{MV}, this is a tensor category equivalent to the category of representations of the Langlands dual group $\check{G}\,\iso\,\GSp_4$. The category $\Sph(\Gr_G)$ acts on the derived category $\D(_{x,\infty}\Bunb_{R_{\pi}})$ by Hecke functors. 
  
  Our main result is Theorem~\ref{Th_main} describing the action of $\Sph(\Gr_G)$ on the irreducible objects of $\P^{\cL}(_{x,\infty}\Bunb_{R_{\pi}})$. It implies the above multiplicity one. It also implies that the action of $\Sph(\Gr_G)$ on $\D(_{x,\infty}\Bunb_{R_{\pi}})$ preserves $\P^{\cL}(_{x,\infty}\Bunb_{R_{\pi}})$. The same phenomenon takes place for Whittaker and Waldspurger models.
 
 To the difference with the case of Whittaker categories, the Bessel category $\P^{\cL}(_{x,\infty}\Bunb_{R_{\pi}})$ is not semi-simple (cf. 2.12). 
  
  The explicit Casselman-Shalika formula for the Bessel models has been established in (\cite{BFF}, Corollary 1.8 and 1.9), where it is presented in the base of $\BM_{\tau}$ consisting of functions supported at a single $R(F)$-orbit on $\Gr_G$. Our Theorem~\ref{Th_main} yields a geometric version of this formula. At the level of functions it yields another base $\{B^{\lambda}\}$ of $\BM_{\tau}$ (cf. 2.14). In this new base the Casselman-Shalika formula writes in an essentially uniform way for Bessel, Waldspurger and Whittaker models.       
     
   In Sect.~1 we propose a general framework that gives a uniform way to define Whittaker, Waldspurger and Bessel categories (the case of Waldspurger models was studied in \cite{Ly2}).

\bigskip\bigskip

{\centerline{\scshape 1. Compactifications and equivariant categories}

\medskip\noindent
1.1  {\scshape Notation} We keep the following notation from \cite{Ly2}. Let $k$ denote an algebraically closed field of characteristic $p\ge 0$. All the schemes (or stacks) we consider are defined over $k$. Let $X$ be a smooth projective connected curve. Fix a
prime $\ell\ne p$. For a scheme (or stack) $S$ write $\D(S)$ for the bounded
derived category of $\ell$-adic \'etale sheaves on $S$, and $\P(S)\subset \D(S)$ for the category of perverse sheaves. 

 Write $\Omega$ for the canonical line bundle on $X$.  For a group scheme $G$ on $X$ write $\cF^0_G$ for the trivial $G$-torsor on $X$.  

\medskip\noindent
1.2.1 Let $G'$ be a connected reductive group over $k$. Given a $G'$-torsor $\gF_{G'}$ on $X$ let $G$ be the group scheme (over $X$) of automorphisms of $\gF_{G'}$. Write $\Bun_G$ for the stack of $G$-bundles on $X$. Note that $\gF_{G'}$ can be viewed as a $G$-torsor as well as a $G'$-torsor on $X$. We identify $\Bun_G$ and $\Bun_{G'}$ via the isomorphism that sends a $G$-torsor $\cF_G$ to the $G'$-torsor $\cF_{G'}=\gF_{G'}\times^G \cF_G$. 

 Let $R\subset G$ be a closed group subscheme over $X$. Say that $G/R$ is \select{strongly quasi-affine over $X$} if for the projection $\pr: G/R\to X$ the $\cO_X$-algebra $\pr_*\cO_{G/R}$ is finitely generated (locally in Zarisky topology), and the natural map $G/R\to \ov{G/R}$ is an open immersion. Here $\ov{G/R}=\Spec (\pr_*\cO_{G/R})$. 

 Let $V$ be a vector bundle on $X$ on which $G$ acts, that is, we are given a homomorphism of group schemes $G\to \Aut(V)$ on $X$. Assume that $R$ is obtained through the following procedure. There is a section $\cO_X\hook{s} V$ such that $V/\cO_X$ is locally free and $R=\{g\in G\mid gs=s\}$. Let $Z$ be the closure of $Gs$ is the total space of $V$, so $G/R \subset Z$. Let $Z'$ be the complement of $Gs$ in $Z$. The following is a consequence of (\cite{Gr}, Theorem~2). 
 
\begin{Lm} 
\label{Lm_qa}
Assume that any fibre of the projection $\pr:Z'\to X$ is of codimenion $\ge 2$ in the corresponding fibre of $\pr:Z\to X$. Then $G/R$ is strongly quasi-affine over $X$, and $Z$ is the affine closure $\ov{G/R}$ of $G/R$. \QED
\end{Lm}
 
 Assume that $R$ satisfies the conditions of Lemma~\ref{Lm_qa} (this holds in our examples below). 
 
\begin{Def} Let $\Bunb_R$ be the following stack. For a scheme $S$ an $S$-point of $\Bunb_R$ is a pair $(\cF_G,\beta)$, where $\cF_G$ is a 
$(S\times X)\times_X G$-torsor on $S\times X$, and $\beta$ is a $G$-equivariant map $\beta: \cF_G\to S\times \ov{G/R}$ over $S\times X$ with the following property. For any geometric point $s\in S$ there is a non-empty open subset $U^s\subset s\times X$ such that 
$$
\beta: \cF_G\mid_{U^s}\to (s\times\ov{G/R})\mid_{U^s}
$$ 
factors through $(s\times G/R)\mid_{U^s}\subset (s\times \ov{G/R})\mid_{U^s}$.
\end{Def}
 
 An $S$-point of $\Bunb_R$ can also be seen as a pair $(\cF_G, \alpha)$, where $\cF_G$ is a $(S\times X)\times_X G$-torsor on $S\times X$, and $\alpha: \cO_{S\times X}\to V_{\cF_G}$ is a section with the following property. First, $\alpha(1)$ lies in $\ov{G/R}\times^G\cF_G$. Secondly, for any geometric point $s\in S$ there is a non-empty open subset $U^s\subset s\times X$ such that $\alpha(1)\mid_{U^s}$ lies in  $(G/R\times^G \cF_G)\mid_{U^s}$. Here $V_{\cF_G}$ is the vector bundle $(V\otimes \cO_{S\times X})\times^G \cF_G$ on $S\times X$. 

  Let $\Bun_R$ denote the stack of $R$-bundles on $X$.
  
\begin{Lm} 
\label{Lm_algebraicity_Bunb_R}
The stack $\Bunb_R$ is algebraic locally of finite type, and 
$\Bun_R\subset \Bunb_R$ is an open substack.
\end{Lm}
\begin{Prf} Consider the stack $\cX$ classifying pairs $(\cF_G,\alpha)$, where $\cF_G$ is a $G$-torsor on $X$, and $\alpha:\cO_X\to V_{\cF_G}$ is a section. It is well-known that this stack is algebraic locally of finite type. The condition that $\alpha(1)$ lies in $\ov{G/R}\times^G\cF_G$ defines a closed substack $\cX'\subset\cX$. The condition that $\alpha(1)$ factors through $G/R\times^G \cF_G$ at the generic point of $X$ is open in $\cX'$.  
Finally, the condition that $\alpha(1)$ lies in $G/R\times^G \cF_G$ everywhere over $X$ is also open.
\end{Prf} 
 
\medskip\noindent
1.2.2  Fix a closed point $x\in X$. Write $\cO_x$ for the completed local ring of $\cO_X$ at $x$, and $F_x$ for its fractions field. 
 
  Let $_{x,\infty}\Bunb_R$ be the following stack. Its $S$-point is a pair $(\cF_G,\alpha)$, where $\cF_G$ is a $(S\times X)\times_X G$-torsor on $S\times X$, and 
\begin{equation}
\label{map_alpha_Bunb_R}
\alpha: \cO_{S\times X}\to V_{\cF_G}(\infty x)
\end{equation}
is a section with the following property. First, $\alpha(1)\mid_{S\times (X-x)}$ lies in 
$
\ov{G/R}\times^G\cF_G\mid_{S\times (X-x)}
$. 
Secondly, for any geometric point $s\in S$ there is a non-empty open subset $U^s\subset s\times (X-x)$ such that 
$\alpha(1)\mid_{U^s}$ lies in  $(G/R\times^G \cF_G)\mid_{U^s}$. 

 Let $\cY_i\subset {_{x,\infty}\Bunb_R}$ be the closed substack given by the condition that (\ref{map_alpha_Bunb_R}) factors through $V_{\cF_G}(ix)\subset V_{\cF_G}(\infty x)$. In particular, $\cY_0=\Bunb_R$.
As in Lemma~\ref{Lm_algebraicity_Bunb_R}, one shows that $\cY_i$ is algebraic locally of finite type. Since $_{x,\infty}\Bunb_R$ is the direct limit of $\cY_i$, the stack $_{x,\infty}\Bunb_R$ is ind-algebraic. 
 
 Remind that if a stack $\cY$ admits a presentation as a direct limit of algebraic stacks locally of finite type $\cY_i$, then we have the derived category $\D(\cY)$, which is an inductive 2-limit of $\D(\cY_i)$. In particular, any $K\in\D(\cY)$ is the extension by zero from some closed algebraic substack of $\cY$. Similarly for the category $\P(\cY)$ of perverse sheaves on $\cY$ (cf. \cite{G}, A.1-A.2 and \cite{BG}, 0.4.4 for details).  
  
  For a scheme $S$, one can also view an $S$-point of $_{x,\infty}\Bunb_R$ as a pair $(\cF_G,\beta)$, where $\cF_G$ is a 
$(S\times X)\times_X G$-torsor on $S\times X$, and $\beta$ is a $G$-equivariant map $\beta: \cF_G\mid_{S\times(X-x)}\to S\times (\ov{G/R}\mid_{X-x})$ with the following property. For any geometric point $s\in S$ there is a non-empty open subset $U^s\subset s\times (X-x)$ such that 
$$
\beta: \cF_G\mid_{U^s}\to (s\times \ov{G/R})\mid_{U^s}
$$ 
factors through $(s\times G/R)\mid_{U^s}\subset (s\times \ov{G/R})\mid_{U^s}$.
  
  Let $H$ be an abelian group scheme over $X$, and $R\to H$ be a homomorphism of group schemes over $X$. Assume that the stack $\Bun_H$ of $H$-bundles on $X$ is algebraic.
  
  Fix a rank one local system $\cL$ on $\Bun_H$ trivialized at the trivial $H$-torsor $\cF^0_H$. Assume that for the tensor product map
$m: \Bun_H\times\Bun_H\to\Bun_H$ there exists an isomorphism $m^*\cL\,\iso\, \cL\boxtimes \cL$ whose restriction to the $k$-point $(\cF^0_H,\cF^0_H)$ is the identity.
    
\medskip\noindent    
1.2.3  We would like to define a category $\P^{\cL}(_{x,\infty}\Bunb_R)$ of \select{$\cL$-equivariant perverse sheaves on $_{x,\infty}\Bunb_R$}, and similarly for $\Bunb_R$. 

  Let $_X\cY\subset (X-x)\times {_{x,\infty}\Bunb_R}$ be the open substack classifying collections $y\in X-x$, $(\cF_G,\beta)\in  {_{x,\infty}\Bunb_R}$ such that the map $\beta:\cF_G\to \ov{G/R}$ factors through $G/R\subset \ov{G/R}$ in a neighbourhood of $y$. 
 
  Set $D_y=\Spec\cO_y$. By definition, for a point of $_X\cY$ the $G$-torsor $\cF_G\mid_{D_y}$ is equipped with a reduction to a $R$-torsor that we denote $\cF_R$.    
  
  Let $_X\cX$ be the stack classifying: $(y, \cF_G,\beta)\in {_X\cY}$, $(y, \cF'_G,\beta')\in {_X\cY}$ and 
\begin{equation}
\label{map_tau_grassmanians}
\tau: \cF_G\mid_{X-y}\,\iso\, \cF'_G\mid_{X-y}
\end{equation}
such that the diagram commutes
$$
\begin{array}{ccc}
\cF_G\mid_{X-y} & \toup{\beta} & \ov{G/R}\mid_{X-y}\\
\downarrow\lefteqn{\scriptstyle \tau} & \nearrow\lefteqn{\scriptstyle \beta'}\\
\cF'_G\mid_{X-y}
\end{array}
$$

 Let $\pr$ (resp., $\act$) denote the projection $_X\cX\to {_X\cY}$ sending the above collection to $(y, \cF_G,\beta)$ (resp., to $(y, \cF'_G,\beta')$). They provide $_X\cX$ with a structure of a groupoid over $_X\cY$.  
  
  Set $D_y^*=\Spec F_y$. Let $_X\GR_R$ denote the stack classifying $(y\in X-x, \cF_R, \cF'_R,\tau)$, where $\cF_R$ and $\cF'_R$ are $R$-torsors on $D_y$ and 
$$
\tau: \cF_R\mid_{D_y^*}\,\iso\, \cF'_R\mid_{D_y^*}
$$ 
is an isomorphism. 

 We have a map $_X\cX\to {_X\GR_R}$ sending the above collection to $(y, \cF_R,\cF'_R,\tau)$, where $\cF_R$ and $\cF'_R$ are $R$-torsors on $D_y$ obtained from $(\cF_G,\beta)$ and $(\cF'_G,\beta')$ and $\tau$ is the restriction of (\ref{map_tau_grassmanians}). 

 Let $_X\Gr_H$ denote the affine grassmanian of $H$ over $X-x$, namely the ind-scheme classifying $y\in X-x$ and an $H$-torsor on $D_y$ trivialized over $D_y^*$. We have a map $_X\GR_R\to{_X\Gr_H}$ sending $(y, \cF_R, \cF'_R,\tau)$ to $(y, \cF_H,\tau)$, where 
$$
\cF_H=\Isom(\cF_R\times_R H, \cF'_R\times_R H),
$$ 
and $\tau:\cF_H\,\iso\, \cF^0_H\mid_{D_y^*}$ is the induced trivialization.  
  
  We have a map $_X\Gr_H\to\Bun_H$ sending $(y, \cF_H,\tau)$ to $\tilde\cF_H$, where $\tilde\cF_H$ is the gluing of $\cF_H^0\mid_{X-y}$ and $\cF_H\mid_{D_y}$ via the isomorphism $\tau: \cF_H\,\iso\,\cF^0_H\mid_{D_y^*}$. 
  
  Define the evaluation map $\ev_{\cX}: {_X\cX}\to \Bun_H$ as the composition
$$
_X\cX\to{_X\GR_R}\to {_X\Gr_H}\to\Bun_H
$$
 
  We would like $\P^{\cL}(_{x,\infty}\Bunb_R)$ to be the category of perverse sheaves $K$ on $_{x,\infty}\Bunb_R$ equipped with an isomorphism
$$
\act^*\tilde K\,\iso\, \pr^*\tilde K\otimes \ev_{\cX}^*\cL
$$
satisfying the usual associativity condition, and such that its restriction to the unit section of $_X\cX$ is the identity. Here $\tilde K$ is the restriction of $K$ under $_X\cY\to{_{x,\infty}\Bunb_R}$. However, this naive definition does not apply directly, because $\pr,\act: {_X\cX}\to {_X\cY}$ are not smooth in general. (One more source of difficulties is that the affine grassmanian $\Gr_{R,y}$ may be highly non-reduced, this happens for example for $R$ a torus). 
    
    We remedy the difficulty under an additional assumption satisfied in our examples. Suppose that $R$ fits into an exact sequence of group schemes $1\to U\to R\to T\to 1$ over $X$, where $U$ is a unipotent group scheme, and $T$ is as follows. There is an integer  $b\ge 0$ and a (ramified) Galois covering $\pi:\tilde X\to X$, where $\tilde X$ is a smooth projective curve,  such that for a $X$-scheme $S$ we have 
$$
T(S)=\Hom(\tilde X\times_X S, \Gm^b)
$$ 
In this case $\Bun_T$ is nothing but the stack of $\Gm^b$-torsors on $\tilde X$. 
For a divisor $D$ on $\tilde X$ with values in the coweight lattice of $\Gm^b$, and for a $T$-torsor $\cF_T$ on $X$, we denote by $\cF_T(D)$ the corresponding twisted $T$-torsor on $X$.

  The stack $_X\cX$ can be seen as the one classifying: $(y, \cF_G, \beta)\in {_X\cY}$, an $R$-torsor $\cF'_R$ on $D_y$, and an isomorphism $\tau: \cF_R\mid_{D^*_y}\,\iso\, \cF'_R\mid_{D^*_y}$, where $\cF_R$ is the $R$-torsor on $D_y$ obtained from $(\cF_G,\beta)$. 
From this point of view the projection $\pr: {_X\cX}\to{_X\cY}$ is the map forgetting $\cF'_R$. 

 Modify the definition of $_X\cX$ and of $_X\cY$ as follows. Let 
$$
_{\tilde X}\cY\subset \tilde X\times {_{x,\infty}\Bunb_R}
$$ 
be the open substack classifying: $\tilde y\in \tilde X$ with $\pi$ nonramified at $\tilde y$ and $y:=\pi(\tilde y)\ne x$, $(\cF_G,\beta)\in  {_{x,\infty}\Bunb_R}$ such that the map $\beta:\cF_G\to \ov{G/R}$ factors through $G/R\subset \ov{G/R}$ in a neighbourhood of $y$. 
 
 Given for each $\sigma\in\Sigma=\Gal(\tilde X /X)$ a coweight $\gamma_{\sigma}:\Gm\to\Gm^b$, we set $\gamma=\{\gamma_{\sigma}\}$. Let  
$$
\pr:{_{\tilde X}\cX_{\gamma}}\to{_{\tilde X}\cY}
$$ 
be the stack whose fibre over $(\tilde y, \cF_G,\beta)\in {_{\tilde X}\cY}$  is the ind-scheme classifying: a $R$-torsor $\cF'_R$ on $D_y$, an isomorphism $\cF_R\,\iso\,\cF'_R\mid_{D^*_y}$, and an extension of the induced isomorphism 
$$
\cF_R\times_R T\,\iso\,\cF'_R\times_R T\mid_{D^*_y}   
$$
to an isomorphism over $D_y$
$$
\cF_R\times_R T\,\iso\,(\cF'_R\times_R T)(\sum_{\sigma \in\Sigma} \gamma_{\sigma} \sigma(\tilde y))
$$  
Here $y=\pi(\tilde y)$, and $\cF_R$ is the $R$-torsor on $D_y$ obtained from $(\cF_G,\beta)$.   
   
  As above, we have an action map $\act:{_{\tilde X}\cX_{\gamma}}\to{_{\tilde X}\cY}$. The advantage is that any fibre of each of the map $\pr,\act: {_{\tilde X}\cX_{\gamma}}\to{_{\tilde X}\cY}$ is reduced (it identifies with the affine grassmanian at $y$ of a unipotent group scheme over $X$). 
  
 Now proceed as in \cite{FGV}. Remind that $U(F_y)$ is an ind-group scheme, it can be written as a direct limit of some group schemes $U^{-m}$, $m\ge 0$, such that $U^{-m}\hook{} U^{-m-1}$ is a closed subgroup, $U^0=U(\cO_y)$, and $U^{-m}/U^0$ are smooth of finite type (\select{loc.cit}, 3.1). 
  
  For this reason, for $m\ge 0$ there exist closed substacks 
$$
_{\tilde X}\cX_{\gamma, m}\hook{} _{\tilde X}\cX_{\gamma, m+1}\hook{} \ldots\hook{}{_{\tilde X}\cX_{\gamma}}
$$ 
such that both maps $\pr, \act: {_{\tilde X}\cX_{\gamma, m}}\to{_{\tilde X}\cY}$ are of finite type and smooth of the same relative dimension, and $_{\tilde X}\cX_{\gamma}$ is a direct limit of the stacks $_{\tilde X}\cX_{\gamma, m}$. 
  
   As above, we have a map $_{\tilde X}\cX_{\gamma}\to {_X\GR_R}$, hence also the evaluation map $\ev_{\tilde X,\gamma}: {_{\tilde X}\cX_{\gamma}}\to \Bun_H$. 
  
\begin{Def} 
\label{Def_P_main}
Let $\P^{\cL}(_{x,\infty}\Bunb_R)$ denote the category of perverse sheaves on $_{x,\infty}\Bunb_R$ equipped for each $\gamma$ and $m\ge 0$ with isomorphisms 
$$
\alpha_{\gamma,m}: \act^*\tilde K\,\iso\, \pr^* \tilde K\otimes\ev^*_{\cX,\gamma}\cL
$$ 
over $_{\tilde X}\cX_{\gamma,m}$. Here $\tilde K$ denotes the restriction of $K$ under $_{\tilde X}\cY\to {_{x,\infty}\Bunb_R}$. It is required that for $m_1<m_2$ the restriction of $\alpha_{\gamma,m_2}$ to $_{\tilde X}\cX_{\gamma, m_1}$ equals $\alpha_{\gamma, m_1}$, the restriction of $\alpha_{0,m}$ to the unit section of $_{\tilde X}\cX_{0,m}$ is the identity, and the usual associativity condition holds. 

 Denote by  $\P^{\cL}(\Bunb_R)$ the full subcategory of $\P^{\cL}(_{x,\infty}\Bunb_R)$ consisting of perverse sheaves, which are extension by zero under $\Bunb_R\hook{}{_{x,\infty}\Bunb_R}$. 
\end{Def}

\medskip\noindent    
1.3   Let $_x\cH_G$ denote the Hecke stack classifying $G$-torsors $\cF_G,\cF'_G$ on $X$ together with an isomorphism $\tau: \cF_G\,\iso\, \cF'_G\mid_{X-x}$. Let $\gq:{_x\cH_G}\to\Bun_G$ (resp., $\gp: {_x\cH_G}\to\Bun_G$) denote the map forgetting $\cF_G$ (resp., $\cF'_G$). 
 Consider the diagram
$$
_{x,\infty}\Bunb_R\getsup{\gp_R} {_{x,\infty}\Bunb_R}\times_{\Bun_G}\,{_x\cH_G}\toup{\gq_R} {_{x,\infty}\Bunb_R},
$$ 
where we used $\gp$ to define the fibred product, $\gp_R$ forgets $\cF'_G$, and $\gq_R$ 
sends 
$
(\cF_G, \beta, \cF'_G,\tau)
$ 
to $(\cF'_G, \beta')$, where $\beta'$ is the composition 
$$
\cF'_G\toup{\tau^{-1}} \cF_G\toup{\beta} \ov{G/R}
$$  
 
  In the same way one gets the diagram
$$
_{\tilde X}\cY \; \getsup{\gp_{\cY}}\;  {_{\tilde X}\cY}\times_{\Bun_G}{_x\cH_G} \;\toup{\gq_{\cY}} \; {_{\tilde X}\cY}
$$
The action of the groupoid $_{\tilde X}\cX$ on $_{\tilde X}\cY$ lifts to an action on this diagram (in the sense of \cite{Ly2}, A.1). Namely, for each $\gamma$ we have two diagrams, where the squares are cartesian
$$
\begin{array}{ccc}
_{\tilde X}\cX_{\gamma} & \toup{\pr} & _{\tilde X}\cY\\
\uparrow\lefteqn{\scriptstyle \gp_{\cX}} && \uparrow\lefteqn{\scriptstyle \gp_{\cY}}\\
_{\tilde X}\cX_{\gamma}\times_{\Bun_G} {_x\cH_G} & \toup{\pr} & _{\tilde X}\cY\times_{\Bun_G}\,{_x\cH_G}\\ 
\downarrow\lefteqn{\scriptstyle \gq_{\cX}} && \downarrow\lefteqn{\scriptstyle \gq_{\cY}}\\
_{\tilde X}\cX_{\gamma} & \toup{\pr} & _{\tilde X}\cY
\end{array}
$$
and 
$$
\begin{array}{ccc}
_{\tilde X}\cX_{\gamma} & \toup{\act} & _{\tilde X}\cY\\
\uparrow\lefteqn{\scriptstyle \gp_{\cX}} && \uparrow\lefteqn{\scriptstyle \gp_{\cY}}\\
_{\tilde X}\cX_{\gamma}\times_{\Bun_G} {_x\cH_G} & \toup{\act} & _{\tilde X}\cY\times_{\Bun_G}\,{_x\cH_G}\\ 
\downarrow\lefteqn{\scriptstyle \gq_{\cX}} && \downarrow\lefteqn{\scriptstyle \gq_{\cY}}\\
_{\tilde X}\cX_{\gamma} & \toup{\act} & _{\tilde X}\cY
\end{array}
$$

 Write $\Sph(\Gr_{G',x})$ for the category of $G'(\cO_x)$-equivariant perverse sheaves on the affine grassmanian $\Gr_{G',x}=G'(F_x)/G'(\cO_x)$. This is a tensor category equivalent to the category of representations of the Langlands dual group $\check{G'}$ over $\Qlb$ (\cite{MV}).

  Let $\Bun_G^x$ be the stack classifying a $G$-bundle $\cF_G$ on $X$ with an isomorphism of $G$-torsors $\cF_G\,\iso\, \,\gF_{G'}\mid_{D_x}$. In a way compatible with our identification $\Bun_G\,\iso\,\Bun_{G'}$ one can view $\Bun_G^x$ as the stack classifying a $G'$-torsor $\cF_{G'}$ with a trivialization $\cF_{G'}\,\iso\, \cF^0_{G'}\mid_{D_x}$. So, the projection $\gq: {_x\cH_G}\to\Bun_G$ can be written as a fibration 
$$
\Bun_G^x\times_{G'(\cO_x)}\Gr_{G',x}\to\Bun_G,
$$
Now for $\cA\in\Sph(\Gr_{G',x})$ and $K\in\D(_{x,\infty}\Bunb_R)$ we can form their twisted exterior product
$$
K\tboxtimes \cA\in \D(_{x,\infty}\Bunb_R\times_{\Bun_G}\,{_x\cH_G})
$$
It is normalized so that it is perverse for $K$ perverse and 
$\DD(K\tboxtimes\cA)\,\iso\, \DD(K)\tboxtimes \, \DD(\cA)$.  

Define the Hecke functor $\H(\cA,\cdot): \D(_{x,\infty}\Bunb_R)\to \D(_{x,\infty}\Bunb_R)$ by 
$$
\H(\cA,K)=(\gp_R)_!(K\tboxtimes \cA)
$$ 
These functors are compatible with the tensor structure on $\Sph(\Gr_{G',x})$. Namely, we have canonically 
\begin{equation}
\label{iso_Hecke_operators}
\H(\cA_1, \H(\cA_2, K))\,\iso\, \H(\cA_1\ast\cA_2, K),
\end{equation}
where $\cA_1\ast\cA_2\in\Sph(\Gr_{G',x})$ is the convolution (\cite{FGV}, Sect.~5).

 As in 1.2, one defines the category $\P^{\cL}(_{x,\infty}\Bunb_R\times_{\Bun_G}\,{_x\cH_G})$. If $K\in\P^{\cL}(_{x,\infty}\Bunb_R)$ then 
$$
K\tboxtimes\cA\in \P^{\cL}(_{x,\infty}\Bunb_R\times_{\Bun_G}\,{_x\cH_G}),
$$ 
so the complex $\H(\cA,K)$ inherits a $\cL$-equivariant structure. Each perverse cohomology sheaf of $\H(\cA,K)$ lies in $\P^{\cL}(_{x,\infty}\Bunb_R)$. 

\medskip\noindent
1.4  Let $\Lambda_{\cY}$ be the set of $R(F_x)$-orbits on the affine grassmanian $\Gr_{G,x}=G(F_x)/G(\cO_x)$. We are interested in the situations where $\Lambda_{\cY}$ is \select{descrete}. Write $\Orb_{\mu}\subset \Gr_{G,x}$ for the $R(F_x)$-orbit corresponding to $\mu\in\Lambda_{\cY}$. 

 Let $\cY_{loc}$ be the stack classifying: a $G$-torsor $\cF_G$ on $D_x$, a $R$-torsor $\cF_R$ on $D^*_x$, and a $R$-equivariant map $\cF_R\to\cF_G\mid_{D_x^*}$. Then $\cY_{loc}$ identifies with the stack quotient of $\Gr_{G,x}$ by $R(F_x)$. 

For $\mu\in\Lambda_{\cY}$ let $\cY^{\mu}_{loc}$ (resp., $\cY^{\le \mu}_{loc}$) denote the stack quotient of $\Orb_{\mu}$ (resp., of $\ov{\Orb}_{\mu}$) by $R(F_x)$. (We don't precise for the moment the scheme structure on $\ov{\Orb}_{\mu}$). We have an order on $\Lambda_{\cY}$ given by $\mu'\le\mu$ iff $\Orb_{\mu'}\subset \ov{\Orb}_{\mu}$. 

 We have a map $_{x,\infty}\Bunb_R\to \cY_{loc}$ sending $(\cF_G,\beta)$ to its restriction to $D_x$. For $\mu\in\Lambda_{\cY}$ set
$$
_{x,\mu}\Bunb_R={_{x,\infty}\Bunb_R}\times_{\cY_{loc}} \cY^{\le\mu}_{loc}\hskip 1em\mbox{and}\hskip 1em
_{x,\mu}\Bunt_R={_{x,\infty}\Bunb_R}\times_{\cY_{loc}} \cY^{\mu}_{loc}
$$ 
Let $_{x,\mu}\Bun_R\subset {_{x,\mu}\Bunt_R}$ be the open substack given by the condition that
$$
\beta: \cF_G\mid_{X-x}\to \ov{G/R}\mid_{X-x}
$$ 
factors through $G/R\mid_{X-x}\subset \ov{G/R}\mid_{X-x}$. 

 To summarize, we have a sequence of embeddings,
$$
_{x,\mu}\Bun_R\hook{} {_{x,\mu}\Bunt_R}\hook{} {_{x,\mu}\Bunb_R}\hook{} {_{x,\infty}\Bunb_R},
$$
where the first two arrows are open embeddings, and the last arrow is a closed one.

\medskip\noindent
1.5 The stack $_{x,\mu}\Bun_R$ classifies: a $G$-torsor $\cF_G$ on $X$, a $G$-equivariant map $\beta: \cF_G\to G/R\mid_{X-x}$ such that the restriction of $(\cF_G,\beta)$ to $D_x$ lies in $\cY^{\mu}_{loc}$.  
 Set 
$$
_{\mu}\cX={_{x,\mu}\Bun_R}\times_{\cY_{loc}}{_{x,\mu}\Bun_R},
$$ 
this is a groupoid over ${_{x,\mu}\Bun_R}$ for the two projections $\pr,\act: {_{\mu}\cX}\to {_{x,\mu}\Bun_R}$. 
 
 View $_{\mu}\cX$ as the stack classifying: $R$-torsors $\cF_R, \cF'_R$ on $X-x$ with an isomorphism $\tau:\cF_R\,\iso\,\cF'_R\mid_{D^*_x}$, a $G$-torsor $\cF_G$ on $X$, and a $R$-equivariant map $\cF_R\to \cF_G\mid_{X-x}$, whose restriction to $D_x$ lies in $\cY_{loc}^{\mu}$. The projection $\pr: {_{\mu}\cX}\to {_{x,\mu}\Bun_R}$ forgets $\cF'_R$. 
 
  Let $_{\mu}\ev_{\cX}: {_{\mu}\cX}\to\Bun_H$ be the map sending the above collection to the $H$-torsor $\tilde\cF_H$ on $X$ obtained by the following gluing procedure. Let $\cF_H$ denote the $H$-torsor on $X-x$ of isomorphisms
$$
\Isom(\cF_R\times_R H, \cF'_R\times_R H)
$$
Then $\tilde\cF_H$ is the gluing of $\cF_H$ and of $\cF^0_H\mid_{D_x}$ over  $D^*_x$ via $\tau: \cF_H\,\iso\, \cF^0_H\mid_{D^*_x}$. 

 We say that $\mu\in\Lambda_{\cY}$ is \select{relevant} if there exists a morphism $\ev^{\mu}:{_{x,\mu}\Bun_R}\to\Bun_H$ making the following diagram commutative
\begin{equation}
\label{diag_relevant}
\begin{array}{ccc}
\Bun_H \times{_{x,\mu}\Bun_R} & \toup{\id\times\ev^{\mu}} \Bun_H\times\Bun_H \toup{m} & \Bun_H\\
\uparrow\lefteqn{\scriptstyle _{\mu}\ev_{\cX}\times\pr} && \uparrow\lefteqn{\scriptstyle \ev^{\mu}}\\
_{\mu}\cX & \toup{\act} & _{x,\mu}\Bun_R,
\end{array}
\end{equation}
If such $\ev^{\mu}$ exists, it is unique up to a tensoring by a fixed $H$-torsor on $X$. Write $\Lambda^+_{\cY}$ for the set of relevant $\mu\in\Lambda_{\cY}$. 

 Write $0\in\Lambda_{\cY}$ for the $R(F_x)$-orbit on $\Gr_{G,x}$ passing by $1$. Then $_{x,0}\Bun_R$ is nothing but the stack $\Bun_R$ of $R$-bundles on $X$. The homomorphism $R\to H$ yields a map $\ev^0:{_{x,0}\Bun_R}\to\Bun_H$ such that (\ref{diag_relevant}) commutes, so $0\in\Lambda^+_{\cY}$.  

  For $\mu\in\Lambda^+_{\cY}$ we denote by $\cB^{\mu}$ the Goresky-MacPherson extension of 
$$
(\ev^{\mu})^*\cL\otimes\Qlb[1](\frac{1}{2})^{\otimes \dim{_{x,\mu}\Bun_R}}
$$
under ${_{x,\mu}\Bun_R}\hook{} {_{x,\mu}\Bunb_R}$. By construction, $\cB^{\mu}\in\P^{\cL}(_{x,\infty}\Bunb_R)$. 

\medskip

  The examples of the above situation include Whittaker models, Waldspurger models for $\GL_2$, and Bessel models for $\GSp_4$ (the latter is studied in Sect.~2). 

\medskip\noindent
1.6 {\scshape Whittaker models} \    Let $G'$ be a connected reductive group over $k$, 
$B'\subset G'$ a Borel subgroup, $U'\subset B'$ its unipotent radical. Set $T'=B'/U'$. Assume that $[G',G']$ is simply connected. Let $\cI$ denote the set of vertices of the Dynkin diagram, and $\{\check{\alpha}_i, i\in\cI\}$ the simple roots corresponding to $B'$. Fix a $B'$-torsor $\gF_{B'}$ on $X$ and a conductor for the induced $T'$-torsor $\gF_{T'}$. That is, for each 
$i\in\cI$ we fix an inclusion of coherent sheaves 
$$
\tilde\omega_i: \cL^{\check{\alpha}_i}_{\gF_{T'}}\hook{}\Omega
$$ 

 Write $\gF_{G'}$ for the $G'$-torsor induced from $\gF_{B'}$. Now $G$ is the group scheme of automorphisms of $\gF_{G'}$. Let $R\subset G$ denote the group scheme of automorphisms of $\gF_{B'}$ acting trivially on $\gF_{T'}$. 
 
 To satisfy the assumptions of Lemma~\ref{Lm_qa}, take
$$
V=\mathop{\oplus}\limits_i \HOM(\cL^{\check{\omega}_i}_{\gF_{T'}}, \, 
\cV^{\check{\omega}_i}_{\gF_{G'}}),  
$$
the sum being taken over the set of fundamental weights $\check{\omega}_i$ of $G'$. Here  $\cV^{\check{\lambda}}$ is the Weil $G'$-module corresponding to $\check{\lambda}$. Then $G$ acts on $V$, and $V$ is equipped with a canonical section $\cO_X\hook{} V$. By (\cite{BG}, 1.1.2), $G/R$ is strongly quasi-affine over $X$. 

   The group scheme of automorphisms of $\gF_{B'/[U',U']}$ acting trivially on $\gF_{T'}$ is canonically 
$$
\mathop{\oplus}\limits_{i\in\cI} \cL^{\check{\alpha}_i}_{\gF_{T'}}
$$ 
Set $H=\mathop{\oplus}\limits_{i\in\cI}\Omega$. Define a homomorphism of group schemes $R\to H$ over $X$ as the composition
$$
R\to \mathop{\oplus}\limits_{i\in\cI} \cL^{\check{\alpha}_i}_{\gF_{T'}}\toup{\tilde\omega} H
$$  
The stack $\Bunb_R$ identifies with the one classifying pairs $(\cF_{G'},\kappa)$, where $\cF_{G'}$ is a $G'$-torsor on $X$, and $\kappa$ is a collection of maps
$$
\kappa^{\check{\lambda}}: \cL^{\check{\lambda}}_{\gF_{T'}}\hook{} \cV^{\check{\lambda}}_{\cF_{G'}}
$$ 
for each dominant weight $\check{\lambda}$ of $G'$, satisfying Pl\" ucker relations (\cite{FGV}, 2.2.2). 

 The set $\Lambda_{\cY}$ identifies in this case with the group $\Hom(\Gm, T')$ of coweights of $T'$. 
 
 For $\lambda\in\Lambda_{\cY}$ the stack $_{x,\lambda}\Bunb_R$ classifies: 
a $G'$-torsor $\cF_{G'}$ on $X$, a collection of maps
$$
\kappa^{\check{\lambda}}: \cL^{\check{\lambda}}_{\gF_{T'}}\hook{} \cV^{\check{\lambda}}_{\cF_{G'}}(\<\lambda,\check{\lambda}\>x)
$$ 
for each dominant weight $\check{\lambda}$ of $G'$, satisfying Pl\" ucker relations.
  
  Assume that the base field $k$ is of characteristic $p>0$, fix a nontrivial additive character $\psi: \Fp\to \Qlb^*$. Write $\cL_{\psi}$ for the corresponding Artin-Shreier sheaf on $\A^1_k$. Take $\cL$ to be the restriction of $\cL_{\psi}$ under the map
$$
\Bun_H\to \prod_{i\in\cI}\H^1(X,\Omega)\toup{\rm{sum}} \A^1_k
$$  
  
   The corresponding Whittaker category $\P^{\cL}(_{x,\infty}\Bunb_R)$ has been described by Frenkel, Gaitsgory and Vilonen in \cite{FGV}. 

\medskip\noindent
1.7 {\scshape Waldspurger models} \  The ground field $k$ is of characteristic $p\ne 2$. Let $\pi:\tilde X\to X$ be a two-sheeted covering ramified over some divisor $D_{\pi}$ on $X$, where $\tilde X$ is a smooth projective curve.  
Set $L_{\pi}=\pi_*\cO_{\tilde X}$ and $G'=\GL_2$. View $L_{\pi}$ as a $G'$-torsor $\gF_{G'}$ on $X$. Let $G$ be the group scheme of automorphisms of $\gF_{G'}$. Let $R$ be the group scheme over $X$ such that for a $X$-scheme $S$ we have 
$
R(S)=\Hom(\tilde X\times_X S, \Gm),
$ 
so $R$ is a closed group subscheme of $G$ over $X$. 

 Let $\sigma$ be the nontrivial automorphism of $\tilde X$ over $X$, so $L_{\pi}\,\iso\, \cO\oplus\cE$, where $\cE$ are $\sigma$-anti-invariants in $L_{\pi}$. It is equipped with $\cE^2\,\iso\,\cO_X(-D_{\pi})$. Take $V=\END_0(L_{\pi})\otimes \cE^{-1}$, where $\END_0(L_{\pi})$ stands for the sheaf of traceless endomorphisms of $L_{\pi}$. The group scheme $G$ acts on $V$ via its action on $L_{\pi}$ (the action of $G$ on $\cE$ is trivial). 
 
 We have 
$$
V\,\iso\, \cO(D_{\pi})\oplus \cO\oplus \cE^{-1}
$$ 
Consider the section $\cO\to V$ given by $(-1,1,0)$. The assumptions of Lemma~\ref{Lm_qa} are satisfied. 

 Set $H=R$. The stack $\Bun_H$ classifies line bundles on $\tilde X$. Pick a rank one local system $\tilde E$ on $\tilde X$. Take $\cL$ to be the automorphic local system on $\Bun_H$ corresponding to $\tilde E$. The stack $_{x,\infty}\Bunb_R$ in this case is canonically isomorphic to the stack $\Wald_{\pi}^x$ introduced in (\cite{Ly2}, 8.2). The corresponding Waldspurger category $\P^{\cL}(_{x,\infty}\Bunb_R)$ has been studied in \select{loc.cit}.

\bigskip\bigskip

{\centerline{\scshape 2. Bessel categories}

\medskip\noindent
2.1.1 {\scshape The group $G$.}  From now on $k$ is an algebraically closed field of characteristic $p>2$. We change the notation compared to Sect.~1. From now on $G=\GSp_4$, so $G$ is the quotient of $\Gm\times\Sp_4$ by the diagonally
embedded $\{\pm 1\}$. We realize $G$ as the subgroup of $\GL(k^4)$ preserving up to a
scalar the bilinear form given by the matrix
$$
\left(
\begin{array}{cc}
0 & E_2\\
-E_2 & 0
\end{array}
\right),
$$
where $E_2$ is the unit matrix of $\GL_2$.

 Let $T$ be the  maximal torus of $G$ given by $\{(y_1,\ldots,y_4)\mid y_iy_{2+i} 
\mbox{ does not depend on  } i\}$. Let $\Lambda$ (resp., $\check{\Lambda}$) denote the
coweight (resp., weight) lattice of $T$.  Let $\check{\epsilon}_i\in\check{\Lambda}$ be the character that sends a point of $T$ to $y_i$. 
 We have $\Lambda=\{(a_1,\ldots,a_4)\in\ZZ^4\mid a_i+a_{2+i} \mbox{ does not depend on }
i\}$ and
$$
\check{\Lambda}=\ZZ^4/\{\check{\epsilon}_1+
\check{\epsilon}_3-\check{\epsilon}_2-\check{\epsilon}_4\} 
$$

 Fix the Borel subgroup of $G$ preserving the flag $ke_1\subset ke_1\oplus ke_2$ of isotropic subspaces in the standard representation. The corresponding positive roots are
$$
\{\check{\alpha}_{12}, \check{\beta}_{ij}, \; 1\le i\le j\le 2\}, 
$$
where $\check{\alpha}_{12}=\check{e}_1-\check{e}_2$ and $\check{\beta}_{ij}=\check{e}_i-\check{e}_{2+j}$.  The simple roots are $\check{\alpha}_{12}$ and $\check{\beta}_{22}$. Write $V^{\check{\lambda}}$ for the irreducible representation of $G$ of highest weight $\check{\lambda}$.

 Fix fundamental weights $\check{\omega}_1=(1,0,0,0)$ and $\check{\omega}_2=(1,1,0,0)$ of $G$. So, $V^{\check{\omega}_1}$ is the standard representation of $G$. The orthogonal to the
coroot lattice is $\ZZ\check{\omega}_0$ with $\check{\omega}_0=(1,0,1,0)$. The orthogonal to the root lattice is $\ZZ\omega$ with $\omega=(1,1,1,1)$.  
 
 Let $P\subset G$ be the Siegel parabolic subgroup preserving the lagrangian subspace $ke_1\oplus ke_2\subset k^4$. Write $U$ for the unipotent radical of $P$, set $M=P/U$. 
 
 Let $\check{G}$ (resp., $\check{M}$)  denote the Langlands dual group over $\Qlb$. Write $V^{\lambda}$ (resp., $U^{\lambda}$) for the irreducible representation of $\check{G}$ (resp., of $\check{M}$) with highest weight $\lambda$. 
 
  Let $w_0$ be the longest element of the Weil group of $G$. Write $\Lambda^+$ for the set of dominant coweights of $G$. By $\check{\rho}$ is denoted the half sum of positive roots of $G$. The corresponding objects for $M$ are denoted $\Lambda^+_M, w_0^M, \check{\rho}_M$. 
  
  Set $G_{ad}=G/Z$, where $Z\subset G$ is the center. Set $\check{\nu}_1=\check{\omega}_2-\check{\omega}_0$ and $\check{\nu}_2=2\check{\omega}_1-\check{\omega}_0$. So, $V^{\check{\nu}_1}$ is the standard representation of $G_{ad}$ and $\wedge^2 V^{\check{\nu}_1}\,\iso\, V^{\check{\nu}_2}$. Let $\Lambda_{G_{ad}}$ be the coweights lattice of $G_{ad}$. Write $\Lambda^{pos}_{G_{ad}}$ for the $\ZZ_+$-span of positive coroots in $\Lambda_{G_{ad}}$. 

\medskip\noindent
2.1.2 For $d\ge 0$ write $X^{(d)}$ for the $d$-th symmetric power of $X$, view it as the scheme of effective divisors of degree $d$ on $X$. Let $^{rss}X^{(d)}\subset X^{(d)}$ denote the open subscheme of divisors of the form $x_1+\ldots+x_d$ with $x_i$ pairwise distinct. Write $\Bun_i$ for the stack of rank $i$ vector bundles on $X$. Set
$$
\RCov^d=\Bun_1\times_{\Bun_1} {^{rss}X^{(d)}},
$$
where the map $^{rss}X^{(d)}\to \Bun_1$ sends $D$ to $\cO_X(-D)$, and the map $\Bun_1\to\Bun_1$ takes a line bundle to its tensor square. It is understood that $^{rss}X^{(0)}=\Spec k$ and the point $^{rss}X^{(0)}\to\Bun_1$ is $\cO_X$. 
Then $\RCov^d$ is the stack classifying two-sheeted coverings $\pi:\tilde X\to X$ ramified exactly at $D\in{^{rss}X^{(d)}}$ with $\tilde X$ smooth (\cite{Ly2}, 7.7.2). 

 Fix a character $\psi:\Fp\to\Qlb^*$ and write $\cL_{\psi}$ for the corresponding Artin-Shreier sheaf on $\A^1$. 

\medskip\noindent
2.2.1 Fix a $k$-point of $\RCov^d$ given by $D_{\pi}\in {^{rss}X^{(d)}}$ and $\pi:\tilde X\to X$ ramified exactly at $D_{\pi}$. Let $\sigma$ denote the nontrivial automorphism of $\tilde X$ over $X$, let $\cE$ be the $\sigma$-anti-invariants in $L_{\pi}:=\pi_*\cO_{\tilde X}$. It is equipped with an isomorphism 
$$
\kappa: \cE^{\otimes 2}\,\iso\, \cO(-D_{\pi})
$$

 Remind that $L_{\pi}$ is equipped with a symmetric form $\Sym^2 L_{\pi}\toup{s}\cO$ such that $\div(L_{\pi}^*/L_{\pi})=D_{\pi}$ for the induced map $L_{\pi}\hook{} L_{\pi}^*$ (\cite{Ly2}, Proposition~14). Set $\cM_{\pi}=L_{\pi}\oplus (L_{\pi}^*\otimes\Omega^{-1})$. It is equipped with a  symplectic form 
 $$ 
\wedge^2 \cM_{\pi}\to L_{\pi}\otimes (L_{\pi}^*\otimes\Omega^{-1})\to \Omega^{-1}
$$
Write $\gF_G$ for the $G$-torsor $(\cM_{\pi}, \Omega^{-1})$ on $X$. 
Let $G_{\pi}$ be the group scheme (over $X$) of automorphisms of $\gF_G$.
Write $\cA_{\pi}$ for the line bundle $\Omega^{-1}$ on $X$ equipped with the corresponding action of $G_{\pi}$.  

 Let $P_{\pi}\subset G_{\pi}$ denote the Siegel parabolic subgroup preserving $L_{\pi}$, $U_{\pi}\subset P_{\pi}$ its unipotent radical. Then $U_{\pi}$ is equipped with a homomorphism of group schemes on $X$
$$
\ev_{\pi}: U_{\pi}\,\iso\, \Omega\otimes\Sym^2 L_{\pi}\toup{s} \Omega
$$ 
Denote by $\tilde R_{\pi}\subset P_{\pi}$ the subgroup stabilizing $\ev_{\pi}$, that is, 
$$
\tilde R_{\pi}=\{p\in P_{\pi}\mid \ev_{\pi}(pup^{-1})=\ev_{\pi}(u)\;\,\mbox{for all}\; u\in U_{\pi}\}
$$
Let $\GL(L_{\pi})$ denote the group scheme (over $X$) of automorphisms of the $\cO_X$-module $L_{\pi}$. Let $T_{\pi}$ denote the functor associating to a $X$-scheme $V$  the group $\H^0(\tilde X\times_X V, \cO^*)$. Then $T_{\pi}$ is a group scheme over $X$, a subgroup of $\GL(L_{\pi})$. 

 Write $\Bun_{T_{\pi}}$ for the stack of $T_{\pi}$-bundles on $X$, that is, for a scheme $S$ the $S$-points of $\Bun_{T_{\pi}}$ is the category of $(S\times X)\times_X T_{\pi}$-torsors on $S\times X$. Given a $\Gm$-torsor on $S\times \tilde X$, its direct image under $\id\times\pi: S\times\tilde X\to S\times X$ is a $(S\times X)\times_X T_{\pi}$-torsor. In this way one identifies $\Bun_{T_{\pi}}$ with the Picard stack $\Pic\tilde X$.  

 Let $\alpha: T_{\pi}\to\Gm$ be the character by which $T_{\pi}$ acts on $\det(L_{\pi})$. Fix an inclusion $T_{\pi}\hook{} \tilde R_{\pi}$ by making $t\in T_{\pi}$ act on $L_{\pi}\oplus (L_{\pi}^*\otimes\Omega^{-1})$ as $(t, \alpha(t)(t^*)^{-1})$, where $t^*\in\Aut(L_{\pi}^*)$ is the adjoint operator.   Set $R_{\pi}=T_{\pi}U_{\pi}$, so $R_{\pi}\subset \tilde R_{\pi}$ is a subgroup. Actually, $\tilde R_{\pi}/U_{\pi}$ identifies with the group of those $g\in\GL(L_{\pi})$ for which there exists $\tilde\alpha(g)\in\Gm$ such that the diagram commutes
$$
\begin{array}{ccc}
\Sym^2 L_{\pi} & \toup{s} & \cO\\ 
 \uparrow\lefteqn{\scriptstyle g} &&  \uparrow\lefteqn{\scriptstyle \tilde\alpha(g)}\\
\Sym^2 L_{\pi} & \toup{s} & \cO
\end{array} 
$$
So, $\tilde R_{\pi}/U_{\pi}$ is equipped with a character $\tilde\alpha: \tilde R_{\pi}/U_{\pi}\to\Gm$ whose restriction to $R_{\pi}$ equals $\alpha$. For $g\in \tilde R_{\pi}/U_{\pi}$ the diagram commutes
$$
\begin{array}{ccc}
L_{\pi} & \toup{s} & L_{\pi}^*\\
 \uparrow\lefteqn{\scriptstyle g} &&  \downarrow\lefteqn{\scriptstyle g^*}\\
L_{\pi} &\toup{\tilde\alpha(g) s} & L_{\pi}^*,
\end{array}
$$
so $(\det g)^2=\tilde\alpha(g)^2$. We see that $R_{\pi}$ is the connected component of $\tilde R_{\pi}$ given by the additional condition $\det g=\tilde\alpha(g)$. 

\begin{Lm} 
\label{Lm_qa_Bess}
The conditions of Lemma~\ref{Lm_qa} are satisfied, so $G_{\pi}/R_{\pi}$ is strongly quasi-affine over $X$. 
\end{Lm}
\begin{Prf}
Define a $G_{\pi}$-module $W_{\pi}$ by the exact sequence $0\to W_{\pi}\to \cA_{\pi}^{-1}\otimes\wedge^2\cM_{\pi}\to\cO_X\to 0$ of $\cO_X$-modules. So, $W_{\pi}$ is equipped with a nondegenerate symmetric form $\Sym^2 W_{\pi}\to\cO$, and the center of $G_{\pi}$ acts trivially on $W_{\pi}$.
  
 We have a subbundle $W_{\pi,1}:=\cA_{\pi}^{-1}\otimes\det L_{\pi}\,\iso\, \Omega\otimes\cE$ in $W_{\pi}$. Let $W_{\pi, -1}$ denote the orthogonal complement to $W_{\pi,1}$ in $W_{\pi}$.  Then $W_{\pi,-1}/W_{\pi,1}\,\iso\, \END_0(L_{\pi})$. As in 1.7, we have a subbundle $\cE\hook{} \END_0(L_{\pi})$. It gives rise to a subbundle
$$
\Omega(-D_{\pi})\hook{} W_{\pi,1}\otimes(W_{\pi,-1}/W_{\pi,1})\hook{}\wedge^2 W_{\pi}
$$ 

 Set
$$
V=(\Omega^{-1}\otimes\cE^{-1}\otimes W_{\pi}) \oplus (\Omega^{-1}(D_{\pi})\otimes\wedge^2 W_{\pi}),
$$
with the action of $G_{\pi}$ coming from its action on $W_{\pi}$. We get a subbundle $\cO_X\hook{s} V$, which is the sum of the above two sections. One checks that $R=\{g\in G\mid gs=s\}$, and the pair $(V,s)$ satisfies the assumptions of Lemma~\ref{Lm_qa}.
\end{Prf}

\medskip\noindent
2.2.2 Fix a $k$-point $x\in X$, write $\cO_x$ for the completed local ring of $X$ at $x$ and $F_x$ for its fraction field. Set $D_x=\Spec\cO_x$ and $D_x^*=\Spec F_x$. 
 
 Write $\tilde F_x$ for the \'etale $F_x$-algebra of regular functions on $\tilde X\times_X D^*_x$. If $x\in D_{\pi}$ then $\tilde F_x$ is nonsplit otherwise it splits over $F_x$. Denote by $\tilde\cO_x$ the ring of regular functions on $\tilde X\times_X D_x$. 

 Write $\Gr_{G_{\pi},x}$ for the affine grassmanian $G_{\pi}(F_x)/G_{\pi}(\cO_x)$. This is an ind-scheme over $k$ that can be seen as the moduli scheme of pairs $(\cF_{G_{\pi}},\beta)$, where
$\cF_{G_{\pi}}$ is a $G_{\pi}$-torsor over $D_x$ and $\beta:\cF_{G_{\pi}}\,\iso\,\cF^0_{G_{\pi}}$ is an isomorphism over $D_x^*$. 

 In concrete terms, $\Gr_{G_{\pi},x}$ classifies pairs: $\cO_x$-lattices $\cM\subset \cM_{\pi}\otimes F_x$ and $\cA\subset \Omega^{-1}\otimes F_x$ such that the diagram commutes 
$$
\begin{array}{ccc}
\wedge^2 \cM_{\pi}\otimes F_x  & \to & \Omega^{-1}\otimes F_x\\
\cup && \cup\\
\wedge^2 \cM & \to & \cA
\end{array}
$$
and induces an isomorphism $\cM\,\iso\, \cM^*\otimes\cA$ of $\cO_x$-modules. 
 
\begin{Def} Let $\cY_{loc}$ denote the stack classifying
\begin{itemize}
\item a free $\tilde F_x$-module $\cB$ of rank one, then write $L$ for $\cB$ viewed as $F_x$-module, it is equipped with the non degenerate form $\Sym^2 L\to\cC$, where $\cC=(\cE\otimes F_x)\otimes \det L$ (\cite{Ly2}, Proposition~14);   

\item a $G$-bundle $(\cM,\cA)$ on $\Spec \cO_x$, here $\cM$ is a free $\cO_x$-module of rank 4, $\cA$ is a free $\cO_x$-module of rank 1 with a symplectic form $\wedge^2 \cM\to\cA$ (it induces $\cM\,\iso\, \cM^*\otimes\cA$);

\item an inclusion $L\hook{} \cM\otimes_{\cO_x} F_x$ of $F_x$-vector spaces, whose image is an isotropic subspace;

\item an isomorphism $\Omega\otimes\cA\otimes F_x\,\iso\, \cC$ of $F_x$-vector spaces.
\end{itemize}
\end{Def} 

\smallskip

\begin{Lm} The stack $\cY_{loc}$ identifies with the stack quotient of $\Gr_{G_{\pi},x}$ by $R_{\pi}(F_x)$. 
\end{Lm}
\begin{Prf}
Given a point of $\cY_{loc}$, it defines a $P_{\pi}$-torsor over $\Spec F_x$. Fix a splitting of the corresponding exact sequence $0\to\Sym^2 L\otimes F_x \to ?\to \cA\otimes F_x\to 0$. Fix also a trivialization $\cB\,\iso\, \tilde F_x$. Then our data becomes just a point of $\Gr_{G_{\pi},x}$. Changing the two trivializations above corresponds to the action of $R_{\pi}(F_x)$ on $\Gr_{G_{\pi},x}$. So, $\cY_{loc}$ classfies a $G_{\pi}$-torsor $\cF_{G_{\pi}}$ on $\D_x$ equipped with a $R_{\pi}$-structure over $D_x^*$.
\end{Prf}

\medskip

  The $R_{\pi}(F_x)$-orbits on $\Gr_{G_{\pi},x}$ are described in (\cite{BFF}, Sect.~1). 
Set $\Lambda_{\cB}=\{(a_1,a_2)\in\ZZ^2\mid a_2\ge 0\}$.  

\begin{Lm} The $k$-points of $\cY_{loc}$ are indexed by $\Lambda_{\cB}$.
\end{Lm}
\begin{Prf}
Given a $k$-point of $\cY_{loc}$, set $L_2=\cM\cap L$. We get a $P_{\pi}$-torsor over $D_x$ given by an exact sequence $0\to\Sym^2 L_2\to ?\to \cA\to 0$ of $\cO_x$-modules. There is a unique $a_1\in\ZZ$ such that the isomorphism over $F_x$ extends to an isomorphism $\Omega\otimes\cA\,\iso\, (\cE\otimes\det L_2)(D_{\pi}+a_1x)$ of $\cO_x$-modules. 

  Further, $(L_2,\cB, L\,\iso\, L_2\otimes F_x)$ is a $k$-point of $\Wald_{\pi}^{x,loc}$ given by some $a_2\ge 0$. Namely, if $\cB_{ex}\subset\cB$ is the smallest $\tilde\cO_x$-lattice such that $L_2\subset \cB_{ex}$ then $a_2=\dim(\cB_{ex}/L_2)$ (\cite{Ly2}, 8.1).
\end{Prf}

\bigskip

 We realize $\Lambda_{\cB}$ as a subsemigroup of $\Lambda_{G_{ad}}$ via the map sending $(a_1,a_2)$ to $\lambda\in\Lambda_{G_{ad}}$ given by 
$\<\lambda, \check{\nu}_1\>=a_1$ and $\<\lambda,\check{\nu}_2\>=a_1+a_2$.  Then $\Lambda_{\cB}=\{\lambda\in\Lambda_{G_{ad}}\mid \<\lambda,\check{\alpha}_{12}\>\ge 0\}$. 

 The image of $\alpha_{12}$ in $\Lambda_{G_{ad}}$ is divisible by two. Define the subsemigroup $\Lambda_{\cB}^{pos}\subset \Lambda_{G_{ad}}$ as the $\ZZ_+$-span of $\frac{1}{2}\alpha_{12}, \beta_{22}$. Then
$$
\Lambda_{\cB}^{pos}=\{\lambda\in\Lambda_{G_{ad}}\mid \<\lambda,\check{\nu}_i\>\ge 0\;\,\mbox{for}\;\, i=1,2\}
$$ 
 
 We introduce an order on $\Lambda_{\cB}$ as follows. For $\lambda,\mu\in\Lambda_{\cB}$ write $\lambda\ge\mu$ iff $\lambda-\mu\in\Lambda_{\cB}^{pos}$. The reader should be cautioned that this is \select{not} the order induced from $\Lambda_{G_{ad}}$ (the latter order is never used in this paper).
 
\medskip\smallskip\noindent
2.3.1 The stack $\Bun_{R_{\pi}}$ classifies collections: a line bundle $\cB_{ex}$ on $\tilde X$, for which we set $L_{ex}=\pi_*\cB_{ex}$, and an exact sequence of $\cO_X$-modules
\begin{equation}
\label{sequence_1}
0\to \Sym^2 L_{ex} \to ?\to \Omega^{-1}\otimes\cE^{-1}\otimes\det L_{ex}\to 0
\end{equation}
By (\cite{Ly2}, Proposition~14), $L_{ex}$ is equipped with a symmetric form 
\begin{equation}
\label{map_quadr_form}
\Sym^2 L_{ex}\to \cE^{-1}\otimes\det L_{ex}
\end{equation}
It admits a canonical section $\cE\otimes\det L_{ex}\hook{s} \Sym^2 L_{ex}$. 

 Here is a Pl\" ucker type description of $\Bun_{R_{\pi}}$. It is the stack classifying:
\begin{itemize}
\item a $G$-bundle $(\cM,\cA)$ on $X$, here $\cM\in\Bun_4$, $\cA\in\Bun_1$ with a symplectic form $\wedge^2 \cM\to\cA$, for which we set $W=\Ker(\cA^{-1}\otimes\wedge^2 \cM\to\cO_X)$;  

\item two subbundles 
$$
\begin{array}{l}
\kappa_1: \Omega\otimes\cE\hook{} W\\
\kappa_2: \Omega(-D_{\pi})\hook{} \wedge^2 W
\end{array}
$$
It is required that there is a lagrangian subbundle $L_{ex}\hook{} \cM$, a line bundle $\cB_{ex}$ on $\tilde X$ and an isomorphism $L_{ex}\,\iso\, \pi_*\cB_{ex}$ with the following properties. Let $W_{-1}$ denote the orthogonal complement to $W_1=\cA^{-1}\otimes\det L_{ex}$ in $W$, so that $W_{-1}/W_1\,\iso\, \END_0(L_{ex})$ is equipped with $\cE\hook{s} \END_0(L_{ex})$. Then

\noindent
$\kappa_1$ factors as $\Omega\otimes\cE\,\iso\, W_1\hook{} W$; 

\noindent
$\kappa_2$ factors as $\Omega(-D_{\pi})\hook{s} W_1\otimes W_{-1}/W_1\hook{} \wedge^2 W$. 
\end{itemize}

\medskip\noindent 
2.3.2 As in 1.2, we have the stacks $\Bunb_{R_{\pi}}\hook{} {_{x,\infty}\Bunb_{R_{\pi}}}$. By definition, $_{x,\infty}\Bunb_{R_{\pi}}$ classifies pairs $(\cF_{G_{\pi}},\beta)$, where $\cF_{G_{\pi}}$ is a $G_{\pi}$-torsor on $X$, and $\beta: \cF_{G_{\pi}}\to \ov{G_{\pi}/R_{\pi}}\mid_{X-x}$ is a $G_{\pi}$-equivariant map such that $\beta$ factors through $G_{\pi}/R_{\pi}$ over some non-empty open subset of $X-x$. 

 Here is a Pl\" ucker type description. The stack ${_{x,\infty}\Bunb_{R_{\pi}}}$ classifies: 
\begin{itemize}
\item a $G$-bundle $(\cM,\cA)$ on $X$, here $\cM\in\Bun_4$, $\cA\in\Bun_1$ with a symplectic form $\wedge^2 \cM\to\cA$, for which we set $W=\Ker(\cA^{-1}\otimes\wedge^2 \cM\to\cO_X)$;  

\item nonzero sections
$$
\begin{array}{l}
\kappa_1: \Omega\otimes\cE\hook{} W(\infty x)\\
\kappa_2: \Omega(-D_{\pi})\hook{} \wedge^2 W(\infty x)
\end{array}
$$
It is required that for some nonempty open subset $X^0\subset X-x$ there is a lagrangian subbundle $L\hook{} \cM\mid_{X^0}$, a line bundle $\cB$ on $\pi^{-1}(X^0)$ and an isomorphism $L\,\iso\, \pi_*\cB\mid_{X^0}$ with the following properties. Let $W_{-1}$ denote the orthogonal complement to $W_1=\cA^{-1}\otimes\det L$ in $W\mid_{X^0}$, so $W_{-1}/W_1\,\iso\, \END_0 L$ is equipped with $\cE\hook{s} \END_0 L$. Then

\noindent
$\kappa_1\mid_{X^0}$ factors as $\Omega\otimes\cE\,\iso\, W_1\hook{} W\mid_{X^0}$; 

\noindent
$\kappa_2\mid_{X^0}$ factors as $\Omega(-D_{\pi})\hook{s} W_1\otimes W_{-1}/W_1\hook{} \wedge^2 W\mid_{X^0}$. 
\end{itemize} 

\begin{Def} For $\lambda\in\Lambda_{\cB}$ denote by $_{x,\lambda}\Bunb_{R_{\pi}}\hook{}{_{x,\infty}\Bunb_{R_{\pi}}}$ the closed substack given by the condition that the maps
$$
\begin{array}{l}
\kappa_1: \Omega\otimes\cE(-\<\lambda, \check{\nu}_1\>x)\hook{} W\\
\kappa_2: \Omega(-D_{\pi}-\<\lambda, \check{\nu}_2\>x)\hook{} \wedge^2 W
\end{array}
$$
initially defined over $X-x$ are regular over $X$. 
\end{Def}

 For $\lambda,\mu\in\Lambda_{\cB}$ we have $_{x,\mu}\Bunb_{R_{\pi}}\subset {_{x,\lambda}\Bunb_{R_{\pi}}}$ 
if and only if $\mu\le\lambda$. 
 
  As in 1.4, we have the open substacks 
$$
 _{x,\lambda}\Bun_{R_{\pi}}\subset {_{x,\lambda}\Bunt_{R_{\pi}}}\subset {_{x,\lambda}\Bunb_{R_{\pi}}},
$$
given by requiring that $\kappa_1, \kappa_2$ are maximal everywhere on $X$ (resp., in a neighbourhood of $x$). 

\medskip\noindent
2.4 The following lemma is straightforward.

\begin{Lm} Let $\lambda\in\Lambda_{\cB}$. For any $k$-point of $_{x,\lambda}\Bunb_{R_{\pi}}$ there is a unique divisor $D$ on $X$ with values in $-\Lambda^{pos}_{\cB}$ such that the maps
$$
\begin{array}{l}
\kappa_1: \Omega\otimes\cE(-\<\lambda x+D, \check{\nu}_1\>)\hook{} W\\
\kappa_2: \Omega(-D_{\pi}-\<\lambda x+D, \check{\nu}_2\>)\hook{} \wedge^2 W
\end{array}
$$
are regular and maximal everywhere on $X$, and $D+\lambda x$ is a divisor with values in $\Lambda_{\cB}$.  \QED
\end{Lm}

 Consider a $\Lambda_{\cB}$-valued divisor $D$ on $X$ with $D=\lambda x+ \sum_{y\ne x} \lambda_y y$ such that $\lambda_y\in -\Lambda_{\cB}^{pos}$ for $y\ne x$. Denote by $_D\Bun_{R_{\pi}}\subset {_{x,\lambda}\Bunb_{R_{\pi}}}$ the substack given by the condition that the maps
$$
\begin{array}{l}
\kappa_1: \Omega\otimes\cE(-\<D, \check{\nu}_1\>)\hook{} W\\
\kappa_2: \Omega(-D_{\pi}-\<D, \check{\nu}_2\>)\hook{} \wedge^2 W
\end{array}
$$
are regular and maximal everywhere on $X$. In particular, for $D=\lambda x$ we get $_D\Bun_{R_{\pi}}\,\iso\, {_{x,\lambda}\Bun_{R_{\pi}}}$.  

  Actually, $_D\Bun_{R_{\pi}}$ is the stack classifying: a line bundle $\cB_{ex}$ on $\tilde X$, for which we set $L_{ex}=\pi_*\cB_{ex}$, a modification 
$L_2\subset L_{ex}$ of rank 2 vector bundles on $X$ such that the composition is surjective
$$
\Sym^2 L_2\to\Sym^2 L_{ex}\to \cE^{-1}\otimes\det L_{ex}
$$
and $\div(L_{ex}/L_2)=\<D, \check{\nu}_2-\check{\nu}_1\>$, and an exact sequence of $\cO_X$-modules
\begin{equation}
\label{eq_seq_in_D_Bun_R}
0\to\Sym^2 L_2\to ?\to \cA\to 0,
\end{equation} 
where
$
\cA=(\Omega^{-1}\otimes\cE^{-1}\otimes \det L_2)(\<D,\check{\nu}_1\>)
$. 
We have used here the description of $\Wald_{\pi}^{x,a}$ from (\cite{Ly2}, Sect.~8.2).

\begin{Rem} 
\label{Rem_a_1}
For $a_1\in\ZZ$ denote by $_x^{a_1}\Bunb_{R_{\pi}}\subset {_{x,\infty}\Bunb_{R_{\pi}}}$ the substack given by the condition that the map
$$
\kappa_1: \Omega\otimes\cE(-a_1 x)\hook{} W
$$
is regular and maximal everywhere on $X$. This is the stack classying collections: $L_2\in\Bun_2$, an exact sequence $0\to\Sym^2 L_2\to ?\to \cA\to 0$ on $X$ with
$\cA=(\Omega^{-1}\otimes\cE^{-1}\otimes\det L_2)(a_1x)$, a line bundle $\cB$ on $\pi^{-1}(X-x)$, and an isomorphism $\pi_*\cB\,\iso\, L_2\mid_{X-x}$. We have the projection 
$$
_x^{a_1}\Bunb_{R_{\pi}}\to \Wald_{\pi}^x
$$ 
sending the above point to $(L_2, \cB, \pi_*\cB\,\iso\, L_2\mid_{X-x})$ (cf. \cite{Ly2}, 8.2). 

 For $\lambda=(a_1, a_2)\in\Lambda_{\cB}$ write $_{x,\lambda}^{a_1}\Bunb_{R_{\pi}}$ for the preimage of $\Wald_{\pi}^{x,\le a_2}$ under this map. The preimage of $\Wald_{\pi}^{x, a_2}$ under the same  map identifies with $_{x,\lambda}\Bun_{R_{\pi}}$. Note that
$$
_{x,\lambda}^{a_1}\Bunb_{R_{\pi}}\subset {_{x,\lambda}\Bunb_{R_{\pi}}}
$$
is an open substack. This will be used in 2.12. 
\end{Rem}  

\medskip\noindent
2.5  Set $H=\Omega\times T_{\pi}$. Denote by $\chi_{\pi}: R_{\pi}\to H$ the homomorphism of group schemes over $X$ given by $\chi_{\pi}(tu)=(ev_{\pi}(u),t)$, $t\in T_{\pi}, u\in U_{\pi}$.  
Let 
$$
ev^0: \Bun_{R_{\pi}}\to \AA^1\times\Pic\tilde X
$$ 
be the map 
sending a point of $\Bun_{R_{\pi}}$ to the pair $(\epsilon, \cB_{ex})$, where $\epsilon$ is the class of the push-forward of (\ref{sequence_1}) by (\ref{map_quadr_form}). 
     
  Fix a rank one local system $\tilde E$ on $\tilde X$. Write $A\tilde E$ for the automorphic local system on $\Pic\tilde X$ corresponding to $\tilde E$. For $d\ge 0$ its inverse image under $\tilde X^{(d)}\to\Pic ^d \tilde X$ identifies with the symmetric power $\tilde E^{(d)}$ of $\tilde E$. 
 
  Let $\cL$ denote the restriction of $\cL_{\psi}\boxtimes A\tilde E$ under the natural map $\Bun_H\to \A^1\times\Pic \tilde X$. As in 1.2, our data give rise to \select{the Bessel category} $\P^{\cL}(_{x,\infty}\Bunb_{R_{\pi}})$. 
   
 One checks that $\lambda=(a_1,a_2)\in\Lambda_{\cB}$ is \select{relevant} (in the sense of 1.4) iff $a_1\ge a_2$. Write $\Lambda_{\cB}^+$ for the set of relevant $\lambda\in\Lambda_{\cB}$.  
 
\medskip\noindent 
2.6 Consider a stratum $_D\Bun_{R_{\pi}}$ of $_{x,\infty}\Bunb_{R_{\pi}}$ as in 2.4, so $D$ is a $\Lambda_{\cB}$-valued divisor on $X$. Arguing as in 1.2.3 (with the difference that now $\tilde y\in\tilde X$ satisfies an additional assumption: $\pi(\tilde y)$ does not lie in the support of $D$), one defines the category $\P^{\cL}(_D\Bun_{R_{\pi}})$. 

 We say that the stratum $_D\Bun_{R_{\pi}}$ is \select{relevant} if $\P^{\cL}(_D\Bun_{R_{\pi}})$ contains a nonzero object. 
As in (\cite{FGV}, Lemma~6.2.8). one shows that the stratum $_D\Bun_{R_{\pi}}$ is relevant iff $D=\lambda x$ with $\lambda\in\Lambda^+_{\cB}$. 
 
 For $\lambda\in\Lambda_{\cB}^+$ denote by 
$$
\ev^{\lambda}: {_{x,\lambda}\Bun_{R_{\pi}}}\to \A^1\times\Pic\tilde X
$$ 
the following map. 
Given a point of  $_{x,\lambda}\Bun_{R_{\pi}}$ as in 2.4, $\ev^{\lambda}$ sends it to the pair $(\epsilon, \cB_{ex})$, where $\epsilon$ is the class of the push-forward of (\ref{eq_seq_in_D_Bun_R}) under the map $\Sym^2 L_2\to\cA\otimes\Omega$, obtained from the symmetric form on $L_{ex}$.  

 For $\lambda\in\Lambda_{\cB}^+$ let $\cB^{\lambda}$ be the Goresky-MacPherson extension of 
$$
(\ev^{\lambda})^*(\cL_{\psi}\boxtimes A\tilde E)\otimes\Qlb[1](\frac{1}{2})^{\otimes \dim{_{x,\lambda}\Bun_{R_{\pi}}}}
$$
under $_{x,\lambda}\Bun_{R_{\pi}}\hook{} {_{x,\lambda}\Bunb_{R_{\pi}}}$. 
The irreducible objects of $\P^{\cL}(_{x,\infty}\Bunb_{R_{\pi}})$ are (up to isomorphism) exactly $\cB^{\lambda}, \lambda\in\Lambda_{\cB}^+$.
 
 Let us underline that for $0\in\Lambda_{\cB}^+$ the only relevant stratum of $_{x,0}\Bunb_{R_{\pi}}=\Bunb_{R_{\pi}}$ is $_{x,0}\Bun_{R_{\pi}}$. So, $\cB^0$ is the extension by zero from $_{x,0}\Bun_{R_{\pi}}$.  As in \cite{FGV}, we say that $\cB^0$ is \select{clean} with respect to the open immersion
$_{x,0}\Bun_{R_{\pi}}\hook{} \Bunb_{R_{\pi}}$. The same argument proves the following.

\begin{Lm} For $\lambda\in\Lambda^+_{\cB}$ the $*$-restriction of $\cB^{\lambda}$ to 
$_{x,\lambda}\Bunt_{R_{\pi}}- {_{x,\lambda}\Bun_{R_{\pi}}}$ vanishes. \QED
\end{Lm}

\medskip\noindent
2.7 The natural projection $\Lambda\to\Lambda_{G_{ad}}$ induces a map
$i:\Lambda^+\to \Lambda^+_{\cB}$. Actually, we get an isomorphism of semi-groups 
$$
\Lambda^+/\ZZ\omega\,\iso\,\Lambda^+_{\cB}
$$ 
The map $i$ preserves the order, that is, if $\lambda\le \mu$ for $\lambda,\mu\in\Lambda^+$ then $i(\lambda)\le i(\mu)$. Besides, $i(-w_0(\lambda))=i(\lambda)$. For $\mu\in\Lambda^+_{\cB}$ an easy calculation shows that 
\begin{equation}
\label{eq_dim_Bess}
 \dim{_{x,\mu}\Bun_{R_{\pi}}}=\<\mu, 2\check{\rho}\>+\dim\Bun_{R_{\pi}}
\end{equation}
 
\begin{Rem} 
\label{Rem_bijection}
Let $\lambda\in\Lambda^+$. The map $\lambda'\mapsto i(\lambda')$ provides a bijection between $\{\lambda' \in\Lambda^+\mid \lambda'\le \lambda\}$ and $\{\mu\in\Lambda^+_{\cB}\mid \mu\le i(\lambda); \; i(\lambda)-\mu=0\;\,\mbox{in}\; \pi_1(G_{ad})\}$.  
\end{Rem}

\medskip\noindent
2.8 Remind that $G=\GSp_4$ and for each $\cA\in\Sph(\Gr_{G,x})$ we have the Hecke  functor $\H(\cA, \cdot): \D(_{x,\infty}\Bunb_{R_{\pi}})\to\D(_{x,\infty}\Bunb_{R_{\pi}})$ introduced in 1.3. 

 Here is our main result. 
 
\begin{Th} 
\label{Th_main}
1) Set $\check{\nu}=\frac{1}{2}w_0(\check{\omega}_0-\check{\beta}_{22})$, so $\check{\nu}\in\check{\Lambda}$. For $\lambda\in\Lambda^+$ we have canonically
$$
\H(\cA_{\lambda}, \cB^0)\,\iso\, \left\{
\begin{array}{ll}
\cB^{i(\lambda)}\otimes(\tilde E_{\tilde x})^{\otimes \<\lambda, 2\check{\nu}\>}, 
& \mbox{the nonsplit case}, \pi(\tilde x)=x\\\\
\cB^{i(\lambda)}\otimes(\tilde E_{\tilde x_1}\otimes\tilde E_{\tilde x_2})^{\otimes \<\lambda, \check{\nu}\>}, & \mbox{the split case}, \pi^{-1}(x)=\{\tilde x_1,\tilde x_2\}
\end{array}
\right.
$$
2) For $\omega=(1,1,1,1)\in\Lambda^+$ and $\mu\in\Lambda^+_{\cB}$ we have canonically
$$
\H(\cA_{\omega}, \cB^{\mu})\,\iso\, \left\{
\begin{array}{ll}
\cB^{\mu}\otimes \tilde E_{\tilde x}^{\otimes 2}, & \mbox{the nonsplit case,}\; \pi(\tilde x)=x\\ \\
\cB^{\mu}\otimes \tilde E_{\tilde x_1}\otimes \tilde E_{\tilde x_2}, & \mbox{the split case,}\; 
\pi^{-1}(x)=\{\tilde x_1, \tilde x_2\}
\end{array}
\right.
$$
\end{Th}

\medskip\noindent
2.9 \  Given a $G$-torsor $\cF_G$ over $D_x$, denote by $\Gr_{G,x}(\cF_G)$ the affine grassmanian classifying pairs $(\cF'_G,\beta)$, where $\cF'_G$ is a $G$-torsor over $D_x$ and $\beta:\cF'_G\,\iso\,\cF_G\mid_{D_x^*}$ an isomorphism. 

 For $\lambda\in\Lambda^+$ we have the subschemes (cf. \cite{BG}, 3.2.1)
$$
\Gr^{\lambda}_{G,x}(\cF_G)\subset \Grb^{\lambda}_{G,x}(\cF_G)\subset \Gr_{G,x}(\cF_G)
$$
A point $(\cF'_G,\beta)\in\Gr_{G,x}(\cF_G)$ lies in $\Grb^{\lambda}_{G,x}(\cF_G)$ if for any $G$-module $V$, whose weights are $\le\check{\lambda}$, we have
$$
V_{\cF_G}(-\<\lambda,\check{\lambda}\>x)\subset V_{\cF'_G}
$$

Remind that we identify $\Gr_{G_{\pi},x}$ with the ind-scheme $\Gr_{G,x}(\gF_G)$ classifying pairs $(\cF_G,\tilde\beta)$, where $\cF_G$ is a $G$-torsor on $D_x$ and 
$$
\tilde\beta:  \cF_G\,\iso\, \gF_G\mid_{D_x^*}
$$ 
is an isomorphism of $G$-torsors. A $k$-point $(\cF_G,\tilde\beta)$ of $\Gr_{G_{\pi},x}$ yields an inclusion $\Grb^{\lambda}_{G,x}(\cF_G)\hook{}\Gr_{G_{\pi},x}$ sending $(\cF'_G,\beta)$ to $(\cF'_G, \tilde\beta\comp\beta)$. 
  For $\mu\in\Lambda_{\cB}$ we denote by $S^{\mu}_{R_{\pi}}\subset \Gr_{G_{\pi},x}$ the  $R_{\pi}(F_x)$-orbit on $\Gr_{G_{\pi},x}$ corresponding to $\mu$. 

 As in \cite{FGV} and (\cite{Ly2}, Proposition~17), the following is a key point of our proof of Theorem~\ref{Th_main}. 
 
\begin{Pp} 
\label{Pp_key}
Let $\mu\in\Lambda^+_{\cB}$. Let $(\cF_G,\tilde\beta)$ be a $k$-point of $S^{\mu}_{R_{\pi}}$, where $\cF_G$ is a $G$-torsor on $D_x$ and $\tilde\beta: \cF_G\,\iso\,\cF^{\pi}_G\mid_{D_x^*}$ is an isomorphism of $G$-torsors. For any $\lambda\in\Lambda^+$ the scheme
\begin{equation}
\label{stratum_key}
\Grb_{G,x}^{\lambda}(\cF_G)\cap S^0_{R_{\pi}}
\end{equation}
is empty unless $\mu\le i(\lambda)$ in the sense of the order on $\Lambda^+_{\cB}$. If $\mu\le i(\lambda)$ then 
\begin{equation}
\label{eq_intersection}
\Gr_{G,x}^{\lambda}(\cF_G)\cap S^0_{R_{\pi}}
\end{equation}
is of dimension $\le \<\lambda,\check{\rho}\>-\<\mu,\check{\rho}\>$. The equality holds if and only if there exists $\lambda'\in\Lambda^+$, $\lambda'\le\lambda$ such that $\mu=i(\lambda')$, and in this case the irreducible components of (\ref{eq_intersection}) of maximal dimension form a base of 
$$
\Hom_{\check{M}}(U^{w_0^{M}w_0(\lambda')}, V^{\lambda})
$$ 
If $\mu=i(\lambda)$ then (\ref{eq_intersection}) is a point scheme. 
\end{Pp}

\begin{Rem} Consider the scheme (\ref{eq_intersection}) in the case $\lambda,\lambda'\in\Lambda^+$ with $\lambda'<\lambda$ and $\mu=i(\lambda')$. Our proof of Proposition~\ref{Pp_key} will also show that for such $\lambda$ and $\mu$ in the nonsplit case \select{all} the irreducible components of (\ref{eq_intersection})  are of the same dimension. In the split case (\ref{eq_intersection}) may have irreducible components of different dimensions (for example, this happens for $\lambda=(a,a,0,0)\in\Lambda^+$ and $\mu=0$). 
\end{Rem}

\noindent
2.10 \  For a $P$-torsor $\cF_P$ over $D_x$ let $\cF_G=\cF_P\times_P G$. For a coweight $\nu\in\Lambda^+_M$ denote by $S^{\nu}_P(\cF_P)$ the ind-scheme classifying pairs $(\cF'_P,\beta)$, where $\cF'_P$ is a $P$-torsor on $D_x$ and 
$$
\beta: \cF'_P\,\iso\,\cF_P\mid_{D_x^*}
$$ 
is an isomorphism such that the pair $(\cF'_M,\beta)$ lies in $\Gr_{M,x}^{\nu}(\cF_M)$. Here $\cF_M$ and $\cF'_M$ are the $M$-torsors induced from $\cF_P$ and $\cF'_P$ respectively. For $\lambda\in\Lambda^+$ denote by 
\begin{equation}
\label{map_t_nu}
\gt^{\nu}_P: S^{\nu}_P(\cF_P)\cap \Gr_{G,x}^{\lambda}(\cF_G) \to \Gr^{\nu}_{M,x}(\cF_M)
\end{equation}
the natural projection. Our Proposition~\ref{Pp_key} is based on the following result established in (\cite{BG}, 4.3.3 and 5.3.7). 

\begin{Pp} 
\label{Pp_borrowed_BG}
All the irreducible components of any fibre of $\gt^{\nu}_P$ are of dimension $\<\nu+\lambda,\check{\rho}\>-\<\nu, 2\check{\rho}_M\>$. These components form a base of
$$
\Hom_{\check{M}}(U^{\nu}, V^{\lambda})
$$ 
For $\nu=w_0^Mw_0(\lambda)$ the map (\ref{map_t_nu}) is an isomorphism. \QED
\end{Pp}

\medskip\smallskip
\begin{Prf}\select{of Proposition~\ref{Pp_key}} 

\smallskip\noindent  
Write $\mu=(a_1,a_2)$. The pair $(\cF_G,\tilde\beta)$ is given by $\cO_x$-lattices $\cM\subset \cM_{\pi}\otimes F_x$ and $\cA\subset \Omega^{-1}\otimes F_x$ such that $(\cM,\cA)$ is a $G$-bundle over $\Spec\cO_x$. Note that 
\begin{equation}
\label{eq_scalar_product}
\<\mu,\check{\rho}\>=\frac{1}{2}(3a_1+a_2)
\end{equation}

\medskip\noindent
1) (the nonsplit case). 

\smallskip
\Step 1   Acting by $R_{\pi}(F_x)$, we may assume that $(\cM,\cA)$ has the following standard form $\cM=L_2\oplus (L_2^*\otimes\cA)$, where $\cA=\Omega^{-1}((a_1-a_2)x)\otimes\cO_x$ and $L_2=\cO_x\oplus \cO_x t^{a_2+\frac{1}{2}}\subset \tilde F_x$, here $t\in\cO_x$ is a local parameter (\cite{Ly2}, Sect.~8.1).
 
  Any $k$-point of $S^0_{R_{\pi}}$ is given by a collection $(a\in\ZZ, L'_2\subset \cM',\cA')$, where $\cM'\subset \cM_{\pi}\otimes F_x$ is a $\cO_x$-lattice,  $\cA'=\Omega^{-1}(-ax)\otimes\cO_x$ and $L'_2=\tilde\cO_x(-a\tilde x)=\cM'\cap (L_{\pi}\otimes F_x)$. Here $\pi(\tilde x)=x$ and $L'_2$ is viewed as a $\cO_x$-module, so 
$$
L'_2=t^{\frac{a}{2}}\cO_x\oplus t^{\frac{a+1}{2}}\cO_x
$$  
Set $\cW=\Ker(\wedge^2 \cM\to\cA)$ and $\cW'=\Ker(\wedge^2 \cM'\to\cA')$. 

 The condition that $(\cF'_G,\beta)=(\cM',\cA')$ lies in $\Grb_{G,x}^{\lambda}(\cF_G)$ implies that
$\cA'\,\iso\,\cA(-\<\lambda,\check{\omega}_0\>x)$, hence 
$$
a=\<\lambda,\check{\omega}_0\>-(a_1-a_2)
$$
It also implies that  
\begin{eqnarray}
\label{inclusion_M}
\cM(-\<\lambda,\check{\omega}_1\>x)\subset \cM'\\
\label{inclusion_W}
\cW(-\<\lambda, \check{\omega}_2\>x)\subset \cW'
\end{eqnarray}
The inclusion (\ref{inclusion_M}) fits into a commutative diagram
$$
\begin{array}{ccccccc}
0\to  & L'_2 & \to & \cM' & \to & {L'}_2^*\otimes\cA' & \to 0\\
       & \cup  && \cup && \cup\\
0\to & L_2(-\<\lambda,\check{\omega}_1\>x) &\to & \cM(-\<\lambda,\check{\omega}_1\>x) & \to & 
L_2^*\otimes\cA(-\<\lambda,\check{\omega}_1\>x) & \to 0
\end{array}
$$ 
This yields an inclusion $L_2^*\subset {L'}^*_2(\<\lambda, \check{\omega}_1-\check{\omega}_0\>)$, which implies $\<\lambda, 2\check{\omega}_1-\check{\omega}_0\>\ge a_1+a_2$. Note that
$2\check{\omega}_1-\check{\omega}_0=\check{\beta}_{12}+\check{\alpha}_{12}$. 

 Further, the inclusion (\ref{inclusion_W}) shows that $(\wedge^2 L_2^*)\otimes\cA^2(-\<\lambda,\check{\omega}_2\>x)\subset (\wedge^2 {L'}_2^*)\otimes\cA'^2$, that is,
$$
\<\lambda, \check{\omega}_2-\check{\omega}_0\>\ge a_1
$$
Since $\check{\omega}_2-\check{\omega}_0=\check{\beta}_{12}$, we get $\mu\le i(\lambda)$. 
 
\medskip 
\Step 2   The above $M$-torsor $(L'_2,\cA')$ is in a position $\nu$ with respect to $(L_2,\cA)$, where $\nu\in\Lambda^+_M$ is a dominant coweight for $M$ that we are going to determine. 

 Clearly, $\<\nu-\lambda,\check{\omega}_0\>=0$. Further, $(\wedge^2 L_2)(-\<\nu,\check{\omega}_2\>x)\,\iso\, \wedge^2 L'_2$, so $a_1=\<\nu, \check{\omega}_0-\check{\omega}_2\>$. From $L_2(-\<\nu,\check{\omega}_1\>x)\subset L'_2$ we get
$$
\<\nu,\check{\omega}_1\>=\left\{
\begin{array}{cl}
 \frac{a}{2}, & a\;\mbox{is even}\\ \\
\frac{a+1}{2}, & a\;\mbox{is odd}
\end{array}
\right.
$$
Now (\ref{eq_intersection}) identifies with the fibre of (\ref{map_t_nu})
over $(L'_2,\cA')\in\Gr^{\nu}_{M,x}(\cF_M)$. Here the $M$-torsor $\cF_M$ is given by $(L_2,\cA)$.

 By Remark~\ref{Rem_bijection}, for $a$ even there exists a unique $\lambda'\in\Lambda^+$ with $\lambda'\le \lambda$ such that $\mu=i(\lambda')$. In this case the above formulas imply $\nu=w_0^Mw_0(\lambda')$. 

 If $\mu=i(\lambda)$ then $a=\<\lambda, \check{\omega}_0-\check{\beta}_{22}\>$ is even, because $\check{\omega}_0-\check{\beta}_{22}$ is divisible by 2 in $\check{\Lambda}$. For $\mu=i(\lambda)$ we get $\nu=w_0^Mw_0(\lambda)$. 
 
 Let us show that $\<\mu,\check{\rho}\>+\<\nu, \check{\rho}-2\check{\rho}_M\>\le 0$. Indeed, since $2\check{\omega}_1-\check{\omega}_2=\check{\alpha}_{12}$, we get
 $$
 \<\nu, \check{\alpha}_{12}\>=\left\{
 \begin{array}{ll}
 a_2, & a\;\mbox{is even}\\
 a_2+1, & a\;\mbox{is odd}
 \end{array}
 \right.
 $$
 and $\<\nu, \check{\alpha}_{12}+\check{\beta}_{22}\>=-a_1$. We have $\check{\rho}-2\check{\rho}_M= \check{\alpha}_{12}+\frac{3}{2}\check{\beta}_{22}$ and $\check{\rho}=2\check{\alpha}_{12}+\frac{3}{2}\check{\beta}_{22}$. So,
$$
\<\nu, \check{\rho}-2\check{\rho}_M\>=\left\{
\begin{array}{ll}
\frac{1}{2}(-3a_1-a_2), & a\;\mbox{is even}\\ \\
\frac{1}{2}(-3a_1-a_2-1), & a\;\mbox{is odd}
\end{array}
\right.
$$
The desired inequality follows now from (\ref{eq_scalar_product}), and it is an equality if and only if $a$ is even, that is, $i(\lambda)-\mu$ vanishes in $\pi_1(G_{ad})$. Our assertion follows now from Proposition~\ref{Pp_borrowed_BG}.

\medskip\noindent
2) (the split case).  

\smallskip
\Step 1
Acting by $R_{\pi}(F_x)$, we may assume that $(\cM,\cA)$ has the following standard form $\cM=L_2\oplus L_2^*\otimes\cA$, where 
$$
L_2=\cO_xt^{a_2}e_1\oplus \cO_x(e_1+e_2)
$$ 
and $\cA=\Omega^{-1}((a_1-a_2)x)\otimes\cO_x$. Here $\{e_i\}$ is a base of $\tilde\cO_x$ over $\cO_x$ consisting of isotropic vectors (\cite{Ly2}, Sect.~8.1). 

 Any $k$-point of $S^0_{R_{\pi}}$ is given by a collection $(b_1,b_2\in\ZZ, L'_2\subset \cM', \cA')$, where $\cM'\subset \cM_{\pi}\otimes F_x$ is a $\cO_x$-lattice, $\cA'=\Omega^{-1}(-(b_1+b_2)x)\otimes\cO_x$ and 
$$
L'_2=\tilde\cO_x(-b_1\tilde x_1-b_2\tilde x_2)=\cM'\cap (L_{\pi}\otimes F_x)
$$ 
Here $\pi^{-1}(x)=\{\tilde x_1,\tilde x_2\}$ and $L'_2$ is viewed as a $\cO_x$-module,so
$$
L'_2=\cO_xt^{b_1}e_1\oplus \cO_xt^{b_2}e_2
$$
If $(\cF'_G,\beta)=(\cM',\cA')$ lies in $\Grb^{\lambda}_{G,x}(\cF_G)$ then $\cA' \,\iso\,\cA(-\<\lambda,\check{\omega}_0\>x)$, so
\begin{equation}
\label{eq_b_1+b_2}
b_1+b_2=\<\lambda,\check{\omega}_0\>-a_1+a_2
\end{equation}

 As in the nonsplit case, the inclusion $L'_2(-\<\lambda, \check{\omega}_1-\check{\omega}_0\>x)\subset L_2$ yields 
\begin{equation}
\label{eq_b_i} 
 b_i+\<\lambda, \check{\omega}_1-\check{\omega}_0\>\ge a_2
\end{equation} 
for $i=1,2$. This implies
$\<\lambda, 2\check{\omega}_1-\check{\omega}_0\>\ge a_1+a_2$. 
 As in the nonsplit case, $(\wedge^2 L'_2)(\<\lambda, 2\check{\omega}_0-\check{\omega}_2\>x)\subset \wedge^2 L_2$ implies $\<\lambda, \check{\omega}_2-\check{\omega}_0\>\ge a_1$. We have shown that $\mu\le i(\lambda)$. 

\medskip
\Step 2  Let us determine $\nu\in\Lambda^+_M$ such that $(L'_2,\cA')\in \Gr^{\nu}_{M,x}(\cF_M)$. Here $\cF_M$ is given by $(L_2,\cA)$.  

 As in the nonsplit case, $\<\nu-\lambda,\check{\omega}_0\>=0$ and $(\wedge^2 L_2)(-\<\nu,\check{\omega}_2\>x)\,\iso\, \wedge^2 L'_2$. So, $a_1=\<\nu, \check{\omega}_0-\check{\omega}_2\>$. From $L_2(-\<\nu,\check{\omega}_1\>x)\subset L'_2$ we get
 $$
 \<\nu,\check{\omega}_1\>=\max\{b_1,b_2\}
 $$

 In particular, for $\mu=i(\lambda)$ we get from (\ref{eq_b_1+b_2}) and (\ref{eq_b_i})
$$
\left\{
\begin{array}{l} 
b_1+b_2=\<\lambda, \check{\omega}_0-\check{\beta}_{22}\>\\
b_i\ge \<\lambda,\check{\alpha}_{12}-\check{\omega}_1+\check{\omega}_0\>
\end{array}
\right.
$$
But $2(\check{\alpha}_{12}-\check{\omega}_1+\check{\omega}_0)=\check{\omega}_0-\check{\beta}_{22}$, so in this case $b_i=\<\lambda,\check{\alpha}_{12}-\check{\omega}_1+\check{\omega}_0\>$ for $i=1,2$. It easily follows that for $\mu=i(\lambda)$ we get $\nu=w_0^Mw_0(\lambda)$.  

 As in the nonsplit case, it remains to show that $\<\mu,\check{\rho}\>+\<\nu, \check{\rho}-2\check{\rho}_M\>\le 0$. We have $\<\nu, \check{\alpha}_{12}+\check{\beta}_{22}\>=-a_1$ and $\<\nu, \check{\alpha}_{12}\>=2\max\{b_i\}-\<\lambda,\check{\omega}_0\>+a_1$. So,
$$
\<\nu, \check{\rho}-2\check{\rho}_M\>=-2a_1-\max\{b_i\}+\frac{1}{2}\<\lambda,\check{\omega}_0\>
$$
The desired inequality follows now from (\ref{eq_scalar_product}), because $\max\{b_i\}\ge \frac{1}{2}(a_2-a_1+\<\lambda,\check{\omega}_0\>)=\frac{1}{2}(b_1+b_2)$.  It is an equality if and only if $b_1=b_2$, and this implies that $2b_i=\<\lambda,\check{\omega}_0\>-(a_1-a_2)$ is even.   

 If $b_1=b_2$ then, as in the nonsplit case, we get $\<\nu, \check{\alpha}_{12}\>=a_2$, so that $\nu=w_0^Mw_0(\lambda')$ for $\lambda'\in\Lambda^+$ such that $\lambda'\le\lambda$ and $i(\lambda')=\mu$.  
\end{Prf} 
 
\medskip\noindent
\begin{Rem} Write $\check{B}\subset\check{G}$ for the dual Borel subgroup in $\check{G}$. The set of double-cosets $\check{M}\backslash \check{G}/\check{B}$ is finite, that is, $\check{M}\subset \check{G}$ is a Gelfand pair. So, for any character $\nu\in\Lambda$ with $\<\nu,\check{\alpha}_{12}\>=0$ and any $\lambda\in\Lambda^+$ the space $\Hom_{\check{M}}(U^{\nu}, V^{\lambda})$ is at most 1-dimensional (\cite{VK}, Theorem~1). This implies that for $\lambda',\lambda\in\Lambda^+$ with $\lambda'\le\lambda$ and $\<\lambda',\check{\alpha}_{12}\>=0$ for $\mu=i(\lambda')$ the scheme (\ref{eq_intersection}) is irreducible.
\end{Rem} 
 
\begin{Rem} 
\label{Rem_Grb_characterization}
Let $\cF_G$ be a $G$-torsor on $D_x$. For a $k$-point $(\cF'_G,\beta)$ of $\Gr_{G,x}(\cF_G)$ we have $(\cF'_G,\beta)\in\Grb^{\lambda}_{G,x}(\cF_G)$ if and only if 
$$
V^{\check{\omega}_i}_{\cF'_G}\subset V^{\check{\omega}_i}_{\cF_G}(\<\lambda, -w_0(\check{\omega}_i)\>x)
$$
for $i=0,1,2$, and for $i=0$ this is an isomorphism.
\end{Rem}

\bigskip\noindent
2.11 Remind the map $\chi_{\pi}: R_{\pi}\to \Omega\times T_{\pi}$ (cf. 2.5).
Write $\chi_{\pi,x}: R_{\pi}(F_x)\to\A^1\times\Pic\tilde X$ for the composition
$$
R_{\pi}(F_x)\toup{\chi_{\pi}} \Omega(F_x)\times T_{\pi}(F_x)\,\iso\,\Omega(F_x)\times\tilde F^*_x\toup{\Res\times\tau_x}\A^1\times\Pic\tilde X,
$$ 
where $\tau_x$ is the natural map $\tilde F_x^*\to \tilde F_x^*/\tilde\cO_x^*\to \Pic\tilde X$.
It is easy to see that for $\mu\in\Lambda^+_{\cB}$ there exists 
a $(R_{\pi}(F_x), \chi_{\pi,x})$-equivariant map $\chi^{\mu}: S^{\mu}_{R_{\pi}}\to\A^1\times\Pic\tilde X$, and such map is unique up to an additive constant (with respect to the structure of an abelian group on $\A^1\times\Pic\tilde X$).

 We need the following analog of (\cite{FGV}, 7.1.7).

\begin{Lm} 
\label{Lm_dominant}
Let $\lambda, \lambda'\in\Lambda^+$ with $\lambda'<\lambda$. Set $\mu=i(\lambda')$. Let $(\cF_G,\tilde\beta)$ be a $k$-point of $S^{\mu}_{R_{\pi}}$. Let $\chi^0: S^0_{R_{\pi}}\to\A^1\times\Pic\tilde X$ be a $(R_{\pi}(F_x), \chi_{\pi,x})$-equivariant map. Then the composition 
\begin{equation}
\label{eq_map_dominant}
\Gr^{\lambda}_{G,x}(\cF_G)\cap S^0_{R_{\pi}}\,\toup{\chi^0}\,\A^1\times\Pic\tilde X\toup{\pr_1}\A^1
\end{equation}
maps each irreducible component of (\ref{eq_intersection}) of dimension $\<\lambda,\check{\rho}\>-\<\mu,\check{\rho}\>$ dominantly to $\A^1$.
\end{Lm}
\begin{Prf}
We may assume that $(\cF_G,\tilde\beta)$ is given by the pair $(\cM,\cA)$ in its standard form as in the proof of Proposition~\ref{Pp_key}, in particular, it is reduced to a $M$-torsor. Write $\mu=(a_1,a_2)$. Set $\nu=w_0^Mw_0(\lambda')$.

Let $z\in\Gm$ act on $L_{\pi}$ as a multiplication by $z$ and trivially on $\Omega^{-1}$. The corresponding action of $\Gm$ on $\cM_{\pi}=L_{\pi}\oplus L_{\pi}^*\otimes\Omega^{-1}$ defines a map $\Gm\to G_{\pi}$ whose image lies in the center of $P_{\pi}/U_{\pi}$. 
 The corresponding action of $\Gm(\cO_x)=\cO_x^*$ on $\Gr_{G_{\pi},x}$ fixes $(\cF_G,\tilde\beta)$ and preserves the scheme (\ref{eq_intersection}). 
  
 The dimension estimates in Proposition~\ref{Pp_key} also show that the irreducible components of dimension $\<\lambda,\check{\rho}\>-\<\mu,\check{\rho}\>$ of the schemes 
$\Gr_{G,x}^{\lambda}(\cF_G)\cap S^0_{R_{\pi}}$ and $\Grb_{G,x}^{\lambda}(\cF_G)\cap S^0_{R_{\pi}}$ are the same. We are going to describe the latter scheme explicitely. 

\medskip\noindent
1) (the split case) 

\medskip\noindent  
We have $\cM=L_2\oplus L_2^*\otimes\cA$ with $L_2=\cO_xt^{a_2}e_1\oplus\cO_x(e_1+e_2)$ and $\cA=\Omega^{-1}((a_1-a_2)x)\otimes\cO_x$, where $\{e_i\}$ is a base of $\tilde\cO_x$ over $\cO_x$ consisting of isotropic vectors, and $t\in\cO_x$ is a local parameter. Let $\cF_M$ be the $M$-torsor on $\Spec\cO_x$ given by $(L_2,\cA)$. 

  Set $b=\frac{1}{2}(a_2-a_1+\<\lambda,\check{\omega}_0\>)$. Consider the $k$-point of $\Gr_{M,x}(\cF_M)$ given by $(L'_2,\cA')$ with  $\cA'=\Omega^{-1}(-2bx)\otimes\cO_x$ and $L'_2=\tilde\cO_x(-b\tilde x_1-b\tilde x_2)$, where $\pi^{-1}(x)=\{\tilde x_1,\tilde x_2\}$. Under our assumptions the scheme (\ref{stratum_key}) identifies with the fibre, say $Y$, of 
\begin{equation}
\label{eq_map_tnu_again}
\gt^{\nu}_P: S^{\nu}(\cF_P)\cap\Grb^{\lambda}_{G,x}(\cF_G)\to \Gr^{\nu}_{M,x}(\cF_M)
\end{equation}
over $(L'_2,\cA')$. In matrix terms, $Y$ is the scheme of those $u\in\Gr_{U,x}$ for which
$gu\in \Grb^{\lambda}_{G,x}$. Here
$$
g=\left(
\begin{array}{cccc}
t^{b-a_2} & -t^{b-a_2} & 0 & 0\\ 
 0 & t^b & 0 & 0\\
 0 & 0 & t^{a_1+b} & 0\\
0 & 0 & t^{a_1-a_2+b} & t^{a_1-a_2+b}
\end{array}
\right)
$$
Write
\begin{equation}
\label{eq_matrix_u}
u=\left(
\begin{array}{cccc}
1 & 0 & u_1 & u_2\\
0 & 1 & u_2 & u_3\\
0 & 0 & 1 & 0\\
0 & 0 & 0 & 1
\end{array}
\right)
\end{equation}
with $u_i\in \Omega(F_x)/\Omega(\cO_x)$. By Remark~\ref{Rem_Grb_characterization}, $Y$ inside of $\Gr_{U,x}$ is given by the equations
$$
\left\{
\begin{array}{l}
u_i\in t^{-b+\<\lambda, w_0(\check{\omega}_1)\>}\Omega(\cO_x)\\\\
 u_i-u_j\in t^{\alpha}\Omega(\cO_x)\\\\
u_1u_3-u_2^2\in  t^{\delta}\Omega^{\otimes 2}(\cO_x)\\\\
u_i\in t^{\delta}\Omega(\cO_x),
\end{array}
\right.
$$
where we have set  for brevity $\delta=-2b+a_2+\<\lambda, w_0(\check{\omega}_2)\>$ and $\alpha=-b+a_2+\<\lambda, w_0(\check{\omega}_1)\>$. 

We may assume that (\ref{eq_map_dominant}) sends (\ref{eq_matrix_u}) to $\Res u_2$. 
Let $Y'\subset Y$ be the closed subscheme given by $u_2=0$. The above $\cO_x^*$-action on $Y$ multiplies each $u_i$ in (\ref{eq_matrix_u}) by the same scalar. So, it sufiices to show that $\dim Y' < \<\lambda,\check{\rho}\>-\<\mu,\check{\rho}\>$. 

 The scheme $Y'$ is contained in the scheme of pairs
$$
\{u_1,u_3\in t^{\delta}\Omega(\cO_x)/\Omega(\cO_x) \mid  u_1u_3\in t^{\delta}\Omega(\cO_x)/\Omega(\cO_x)\}
$$
The dimension of the latter scheme is at most $-\delta$. We have $-\delta\le  \<\lambda,\check{\rho}\>-\<\mu,\check{\rho}\>$, and the equality holds iff $\alpha=0$. But if $\alpha=0$ then $Y'$ is a point scheme. Since $\<\lambda,\check{\rho}\>-\<\mu,\check{\rho}\>$ is strictly positive, we are done.

\medskip
\noindent
2) (the nonsplit case) 
 
\medskip\noindent   
We have $\cM=L_2\oplus L_2^*\otimes\cA$ with $L_2\,\iso\,\cO_x\oplus \cO_xt^{a_2+\frac{1}{2}}$ and $\cA\,\iso\,\Omega^{-1}((a_1-a_2)x))\otimes\cO_x$, where $t\in\cO_x$ is a local parameter. 
Let $\cF_M$ be the $M$-torsor on $\Spec\cO_x$ given by $(L_2,\cA)$.

  Set $L_2'=t^{\frac{a}{2}}\cO_x\oplus t^{\frac{a+1}{2}}\cO_x$ and $\cA'=\Omega^{-1}(-ax)\otimes\cO_x$ with $a=\<\lambda,\check{\omega}_0\>-a_1+a_2$,  remind that $a$ is even.  
The scheme (\ref{stratum_key}) identifies with the fibre, say $Y$, of (\ref{eq_map_tnu_again})
over $(L'_2,\cA')$.
 
  Consider the base $\{1,t^{\frac{1}{2}}\}$ in $L_{\pi}\otimes\cO_x$ and the dual base in $L_{\pi}^*\otimes\cO_x$. Then in matrix terms, $Y$ becomes the scheme of those $u\in\Gr_{U,x}$ for which $gu\in \Grb^{\lambda}_{G,x}$. Here
$g=t^{\frac{a}{2}}\diag(1, t^{-a_2}, t^{a_1-a_2}, t^{a_1})$. For $u\in\Gr_{U,x}$ written as (\ref{eq_matrix_u}), the scheme $Y$ is given by the equations
$$
\left\{
\begin{array}{l}
u_1\in t^{-\frac{a}{2}+\<\lambda, w_0(\check{\omega}_1)\>}\Omega(\cO_x)\\\\
u_2,u_3\in  t^{\alpha}\Omega(\cO_x)\\\\
u_1u_3-u_2^2\in t^{\delta}\Omega^{\otimes 2}(\cO_x)\\\\
u_i\in t^{\delta}\Omega(\cO_x),
\end{array}
\right.
$$
where we have set $\alpha=a_2-\frac{a}{2}+\<\lambda, w_0(\check{\omega}_1)\>$ and $\delta=a_2-a+\<\lambda, w_0(\check{\omega}_2)\>$.

We may assume that   (\ref{eq_map_dominant}) sends (\ref{eq_matrix_u}) to $\Res (u_1-tu_3)$. 
Let $Y'\subset Y$ be the closed subscheme given by $u_1=tu_3$. Since we have an action of $\cO_x^*$, it suffices to show that $\dim Y'< \<\lambda,\check{\rho}\>-\<\mu,\check{\rho}\>$. 
 
  The scheme $Y'$ is contained in the scheme 
 $$
\{u_2,u_3\in t^{\delta}\Omega(\cO_x)/\Omega(\cO_x) \mid tu_3^2-u_2^2\in t^{\delta}\Omega^{\otimes 2}(\cO_x)\}
 $$
The latter scheme is included into $Y''$ given by
 $$
 \begin{array}{ll}
Y'' = \{u_2,u_3\in t^{\frac{\delta}{2}}\Omega(\cO_x)/\Omega(\cO_x)\}, & \mbox{for $\delta$ even}\\\\
Y'' = \{u_2\in t^{\frac{1+\delta}{2}}\Omega(\cO_x)/\Omega(\cO_x), \; u_3\in t^{\frac{\delta -1}{2}}\Omega(\cO_x)/\Omega(\cO_x)\}, & \mbox{for $\delta$ odd}
 \end{array}
$$
This implies $\dim Y'\le \dim Y''\le-\delta$. As in the split case, $-\delta\le \<\lambda,\check{\rho}\>-\<\mu,\check{\rho}\>$ and the equality implies $\alpha=0$. But for $\alpha=0$ we get $Y'\,\iso\,\Spec k$. This concludes the proof.
\end{Prf}

\bigskip\noindent
\begin{Prf}\select{of Theorem~\ref{Th_main}}

\medskip\noindent
2) Let $\gq_{\omega}:{_{x,\infty}\Bunb_{R_{\pi}}}\,\iso\, {_{x,\infty}\Bunb_{R_{\pi}}}$ denote the isomorphism sending $(\cM,\cA, \kappa_1,\kappa_2)$ to 
$$
(\cM(x),\cA(2x), \kappa_1,\kappa_2)
$$ 
It preserves the stratification of $_{x,\infty}\Bunb_{R_{\pi}}$ by $_D\Bun_{R_{\pi}}$ introduced in 2.4, and we have a commutative diagram
$$
\begin{array}{ccc}
_{x,\mu}\Bun_{R_{\pi}} & \toup{\gq_{\omega}} & {_{x,\mu}\Bun_{R_{\pi}}} \\
\downarrow\lefteqn{\scriptstyle \ev^{\mu}} && \downarrow\lefteqn{\scriptstyle \ev^{\mu}} \\
\AA^1\times\Pic \tilde X & \toup{\id\times \,\tilde\gq_{\omega}} & \AA^1\times\Pic \tilde X,
\end{array}
$$
where $\tilde\gq_{\omega}$ sends $\cB_{ex}$ to $\cB_{ex}(2\tilde x)$ (resp., to $\cB_{ex}(\tilde x_1+\tilde x_2)$) in the nonsplit (resp., split) case. Our assertion follows from the automorphic property of $A\tilde E$.

\medskip\noindent
1) We change the notation replacing $\lambda$ by $-w_0(\lambda)$. In other words, we will establish a canonical isomorphism $\H(\cA_{-w_0(\lambda)}, \cB^0)\,\iso\, \cB^{i(\lambda)}\otimes\cN$ with
$$
\cN\,\iso\, \left\{
\begin{array}{ll}
(\tilde E_{\tilde x})^{\otimes \<\lambda, 2\check{\nu}\>}, 
& \mbox{the nonsplit case}, \pi(\tilde x)=x\\\\
(\tilde E_{\tilde x_1}\otimes\tilde E_{\tilde x_2})^{\otimes \<\lambda, \check{\nu}\>}, & \mbox{the split case}, \pi^{-1}(x)=\{\tilde x_1, \tilde x_2\}
\end{array}
\right.
$$
Denote by $\tilde K_{\mu}$ (resp., by $K_{\mu}$, $_D K$) the $*$-restriction of $\H(\cA_{-w_0(\lambda)}, \cB^0)$ to $_{x,\mu}\Bunt_{R_{\pi}}$ (resp., to $_{x,\mu}\Bun_{R_{\pi}}$, $_D\Bun_{R_{\pi}}$). Here $D$ is $\Lambda_{\cB}$-valued divisor on $X$ as in 2.4. 

By 1.3, we know that each perverse cohomology sheaf of $_DK$ is $\cL$-equivariant. So, $_DK=0$ unless $D=\mu x$ with $\mu$ relevant. In partiuclar, $\tilde K_{\mu}$ is the extension by zero under $_{x,\mu}\Bun_{R_{\pi}}\hook{} {_{x,\mu}\Bunt_{R_{\pi}}}$. 

 Since $\cB^0$ is self-dual (up to replacing $\tilde E$ by $\tilde E^*$ and $\psi$ by $\psi^{-1}$), our assertion is reduced to the following lemma.
\end{Prf}

\begin{Lm} We have $\tilde K_{\mu}=0$ unless $\mu\le i(\lambda)$. The complex $\tilde K_{\mu}$ lives in non positive (resp., strictly negative) perverse degrees for $\mu=i(\lambda)$ (resp., for $\mu< i(\lambda)$). We have canonically
$$
K_{i(\lambda)}\,\iso\, (\ev^{i(\lambda)})^*(\cL_{\psi}\boxtimes A\tilde E)\otimes\cN\otimes \Qlb[1](\frac{1}{2})^{\otimes \dim{_{x,i(\lambda)}\Bun_{R_{\pi}}}}
$$
\end{Lm}
\begin{Prf}
Write $_x\overl{\cH}^{\lambda}_G$ for the substack of $_x\cH_G$ that under the projection $\gq_G:{_x\cH_G}\to\Bun_G$ identifies with 
$$
\Bun_G^x\times_{G(\cO_x)}\Grb^{\lambda}_{G,x}\to \Bun_G
$$
 
 For the diagram 
$$
_{x,\infty}\Bunb_{R_{\pi}}\,\getsup{\gp_R}\, {_{x,\infty}\Bunb_{R_{\pi}}}\times_{\Bun_G}{_x\overl{\cH}^{-w_0(\lambda)}_G}\,
\toup{\gq_R} \,{_{x,\infty}\Bunb_{R_{\pi}}}
$$
we have 
$$
\H(\cA_{-w_0(\lambda)},\; \cdot)=(\gp_R)_!(\gq_R^*(\cdot)\tboxtimes\cA_{-w_0(\lambda)})
$$ 

 Let $\mu=(a_1,a_2)\in\Lambda^+_{\cB}$. Pick a $k$-point $\eta\in {_{x,\mu}\Bun_{R_{\pi}}}$ given by a collection:  a line bundle $\cB_{ex}$ on $\tilde X$, for which we set $L_{ex}=\pi_*\cB_{ex}$, a modification $L_2\subset L_{ex}$ of rank 2 vector bundles on $X$ such that the composition is surjective
$$
\Sym^2 L_2\to \Sym^2 L_{ex}\to (\cE\otimes\det L_{ex})(D_{\pi})
$$
and $a_2x=\div(L_{ex}/L_2)$, and an exact sequence 
\begin{equation}
\label{sequence_local_P}
0\to\Sym^2 L_2\to ?\to\cA\to 0
\end{equation}
on $X$, where we have set $\cA=(\Omega^{-1}\otimes\cE\otimes\det L_2)(D_{\pi}+a_1x)$.  

  The fibre of 
$$
\gp_R: {_{x,\infty}\Bunb_{R_{\pi}}}\times_{\Bun_G}{_x\overl{\cH}^{-w_0(\lambda)}_G}
\to {_{x,\infty}\Bunb_{R_{\pi}}}
$$ 
over $\eta$ identifies with $\Grb_{G,x}^{\lambda}(\cF_G)$, where $\cF_G=(\cM,\cA)\in\Bun_G$ is given by the $P$-torsor (\ref{sequence_local_P}). 
  
  Fix a trivialization $\cB_{ex}\otimes\tilde\cO_x\,\iso\,\tilde\cO_x$ and a splitting of (\ref{sequence_local_P}) over $\Spec\cO_x$. They yield isomorphisms $\cM\,\iso\, (L_2\oplus L_2^*\otimes\cA)\mid_{\Spec\cO_x}$ and $\cA\,\iso\,\Omega^{-1}((a_1-a_2)x)\mid_{\Spec\cO_x}$. So, the pair 
$$
\left\{
\begin{array}{c}
\cM\otimes\cO_x\subset \cM_{\pi}\otimes F_x\\
\cA\otimes\cO_x\subset \Omega^{-1}\otimes F_x
\end{array}
\right.
$$
becomes a point of $\Gr_{G_{\pi},x}$ lying in $S^{\mu}_{R_{\pi}}$.  

 Remind that $\cB^0$ is clean with respect to the open immersion $_{x,0}\Bun_{R_{\pi}}\subset {_{x,0}\Bun_{R_{\pi}}}$. So, only the stratum (\ref{stratum_key}) contributes to $K_{\mu}$. By Proposition~\ref{Pp_key}, $K_{\mu}=0$ unless $\mu\le i(\lambda)$. 

 Assume that $\mu\le i(\lambda)$. Stratify (\ref{stratum_key}) by locally closed subschemes
$\Gr_{G,x}^{\lambda'}(\cF_G)\cap S^0_{R_{\pi}}$ with $\lambda'\le\lambda$, where $\lambda'\in\Lambda^+$. The $*$-restriction of $\cA_{-w_0(\lambda)}$ under 
$$
\Gr_{G,x}^{\lambda'}(\cF_G)\hook{}\Grb_{G,x}^{\lambda}(\cF_G)
$$ 
is a constant complex placed in usual degree $\le -\dim\Gr_{G,x}^{\lambda'}(\cF_G)=-\<\lambda', 2\check{\rho}\>$, the inequality is strict unless $\lambda'=\lambda$.  From (\ref{eq_dim_Bess}) and Proposition~\ref{Pp_key}, we get
$$
-\dim{_{x,0}\Bun_{R_{\pi}}} -\<\lambda', 2\check{\rho}\>+2\dim(\Gr^{\lambda'}_{G,x}(\cF_G)\cap S^0_{R_{\pi}}) \le -\dim {_{x,\mu}\Bun_{R_{\pi}}}
$$
So, $K_{\mu}$ is placed in perverse degrees $\le 0$. If $\mu-i(\lambda)$ does not vanish in $\pi_1(G_{ad})$ then, by Proposition~\ref{Pp_key}, $K_{\mu}$ is placed is strictly negative perverse degrees.

  If $i(\lambda)-\mu$ vanishes in $\pi_1(G_{ad})$ let $\lambda'\in\Lambda^+$ be such that $\lambda'\le\lambda$ and $\mu=i(\lambda')$. Then only the stratum (\ref{eq_intersection}) could contribute to the 0-th perverse cohomology sheaf of $K_{\mu}$. For $\mu<i(\lambda)$ it does not contribute, because the restriction of $\gq^*_R(\cB^0)\tboxtimes\cA_{-w_0(\lambda)}$ to (\ref{eq_intersection}) is a nonconstant local system by Lemma~\ref{Lm_dominant}.

 If $\mu=i(\lambda)$ then (\ref{eq_intersection}) is a point scheme by Proposition~\ref{Pp_key}, and the description of $K_{i(\lambda)}$ follows from the automorphic property of $A\tilde E$.
\end{Prf}


\medskip\smallskip\noindent
2.12 For $\lambda\in\Lambda^+_{\cB}$ the perverse sheaf $\cB^{\lambda}$ is not always the extension by zero from $_{x,\lambda}\Bun_{R_{\pi}}$. For example, take $\lambda=(1,1)$ and $\mu=(1,0)$. An easy calculation shows that, over $_{x,\lambda}\Bun_{R_{\pi}}\cup {_{x,\mu}\Bun_{R_{\pi}}}$, $\cB^{\lambda}$ is a usual sheaf placed in cohomological degree 
$-\dim {_{x,\lambda}\Bun_{R_{\pi}}}$.

  Now we can show that the category $\P^{\cL}(_{x,\infty}\Bunb_{R_{\pi}})$ is not semi-simple. Remind the stack $_x^{a_1}\Bunb_{R_{\pi}}$ (cf. Remark~\ref{Rem_a_1}). Let $\lambda=(1,1)$ and $\mu=(1,0)$. We have a sequence of open embeddings
$$
{_{x,\lambda}\Bun_{R_{\pi}}}\hook{j}\; {_{x,\lambda}^1\Bunb_{R_{\pi}}} \;\hook{\tilde j}\;
{_{x,\lambda}\Bunb_{R_{\pi}}},
$$ 
where $j$ is obtained from the affine open embedding 
$\Wald^{x,1}_{\pi}\hook{} 
\Wald^{x,\le 1}_{\pi}$ 
by the base change 
$$
{_{x,\lambda}^1\Bunb_{R_{\pi}}}\to \Wald^{x,\le 1}_{\pi}
$$ 
Set $\cB^{\lambda, \mu}=\tilde j_{! *} j_!(\cB^{\lambda}\mid_{_{x,\lambda}\Bun_{R_{\pi}}})
$. 
We get an exact sequence in $\P^{\cL}(_{x,\infty}\Bunb_{R_{\pi}})$
$$
0\to K\to \cB^{\lambda,\mu}\to\cB^{\lambda}\to 0
$$
If $\P^{\cL}(_{x,\infty}\Bunb_{R_{\pi}})$ was semi-simple, it would split, this contradicts the fact that the $*$-restriction of $\cB^{\lambda}$ to $_{x,\mu}\Bun_{R_{\pi}}$ is not zero. 

\bigskip\noindent
2.13 {\scshape Geometric Casselman-Shalika formula} 

\smallskip\noindent
Remind that we write $V^{\mu}$ for the irreducible representation of $\check{G}$ of highest weight $\mu$. Let $E$ be a $\check{G}$-local system on $\Spec k$ equipped with an isomorphism
$$
V^{\omega}_E\,\iso\, \left\{
\begin{array}{ll}
\tilde E_{\tilde x}^{\otimes 2}, & \mbox{the nonsplit case,}\; \pi(\tilde x)=x\\ \\
\tilde E_{\tilde x_1}\otimes \tilde E_{\tilde x_2}, & \mbox{the split case,}\; 
\pi^{-1}(x)=\{\tilde x_1,\tilde x_2\}
\end{array}
\right.
$$
We assign to $E$ the ind-object $K_E$ of $\P^{\cL}(_{x,\infty}\Bunb_{R_{\pi}})$ given by 
$$
K_E=\mathop{\oplus}\limits_{
\begin{array}{c}
\scriptstyle{\lambda\in\Lambda^+}\\
\scriptstyle{\<\lambda,\check{\nu}\>=0}
\end{array}}
\cB^{i(\lambda)}\otimes (V^{\lambda})^*_E,  
$$
where $\check{\nu}\in\check{\Lambda}$ is that of Theorem~\ref{Th_main}. For a representation $V$ of $\check{G}$ write $\cA_V$ for the object of $\Sph(\Gr_{G,x})$ corresponding to $V$ via the Satake equivalence $\Rep(\check{G})\,\iso\,\Sph(\Gr_{G,x})$. 

 One formally derives from Theorem~1 the following.
\begin{Cor} 
For any $V\in\Rep(\check{G})$ there is an isomorphism
$\alpha_V: \H(\cA_V, K_E)\,\iso\, K_E\otimes V_E$. For $V,V'\in\Rep(\check{G})$ the diagram commutes
$$
\begin{array}{ccc}
\H(\cA_{V'}, \H(\cA_V, K_E)) & \toup{\alpha_V} & \H(\cA_{V'}, K_E\otimes V_E)\\ 
\downarrow\lefteqn{\scriptstyle \eta} && \downarrow\lefteqn{\scriptstyle\alpha_{V'}\otimes\id}\\
\H(\cA_{V\otimes V'}, K_E) & \toup{\alpha_{V\otimes V'}} & K_E\otimes (V\otimes V')_E,
\end{array}
$$
where $\eta$ is the isomorphism (\ref{iso_Hecke_operators}).
\end{Cor}

\medskip\noindent
2.14 One may view $\Gr_{G_{\pi},x}$ as the ind-scheme classifying a $G_{\pi}$-bundle $\cF_{G_{\pi}}$ on $X$ together with a trivialization $\cF_{G_{\pi}}\,\iso\, \cF^0_{G_{\pi}}\mid_{X-x}$. This yields a map $\Gr_{G_{\pi},x}\to {_{x,\infty}\Bunb_{R_{\pi}}}$. 

 Theorem~\ref{Th_main} holds also in the case of a finite base field $k=\Fq$. In this case we have the Bessel module $\BM_{\tau}$ introduced in 0.1, which we now view as the space of  functions on $G_{\pi}(F_x)/G_{\pi}(\cO_x)$ that change by $\tau$ under the action of $R_{\pi}(F_x)$. Let $B^{\lambda}$ denote the restriction under 
$$
G_{\pi}(F_x)/G_{\pi}(\cO_x)\to {_{x,\infty}\Bunb_{R_{\pi}}}(k)
$$ 
of the trace of Frobenius function of $\cB^{\lambda}$. Then $\{B^{\lambda}, \lambda\in\Lambda^+_{\cB}\}$ is a base of $\BM_{\tau}$. From Theorem~1 it follows that $\BM_{\tau}$ is a free module of rank one over the Hecke algebra $\H_{\chi_c}$.

\end{document}